\documentclass[final]{siamart1116}
\usepackage{amssymb}
\usepackage{eurosym}
\usepackage{lipsum}
\usepackage{amsfonts}
\usepackage{graphicx}
\usepackage{algorithmic}
\usepackage{amsopn}

\ifpdf
  \DeclareGraphicsExtensions{.eps,.pdf,.png,.jpg}
\else
  \DeclareGraphicsExtensions{.eps}
\fi
\numberwithin{theorem}{section}
\newcommand{\TheTitle}{Lavrent'ev equation}
\newcommand{\TheAuthors}{M.V.Klibanov J.Li and W.Zhang}
\newtheorem{example}{Example}
\headers{\TheTitle}{\TheAuthors}

\ifpdf
\hypersetup{
  pdftitle={\TheTitle},
  pdfauthor={\TheAuthors}
}
\fi

\renewcommand{\thealgorithm}{\arabic{section}.\arabic{algorithm}}
\renewcommand{\thefigure}{\arabic{section}.\arabic{figure}}
\renewcommand\theequation{\thesection.\arabic{equation}}
\input{tcilatex}
\begin{document}

\title{Linear Lavrent'ev Integral Equation for the Numerical Solution of a
Nonlinear Coefficient Inverse Problem\thanks{%
Submitted to the editors DATE. 
\funding{The effort of Klibanov was supported by US Army
Research Laboratory and US Army Research Office grant W911NF-19-1-0044. The
work of Li was partially supported by the NSF of China No. 11971221 and the
Shenzhen Sci-Tech Fund No. JCYJ20170818153840322 and JCYJ20190809150413261.
The work of Zhang was partially supported by the Shenzhen Sci-Tech Fund No.
JCYJ20170818153840322, JCYJ20180307151603959 and the NSF of China No.
11901282 and 11731006.}}}
\author{ Michael V. Klibanov\thanks{%
Department of Mathematics and Statistics, University of North Carolina at
Charlotte, Charlotte, NC 28223, USA (mklibanv@uncc.edu)} \and Jingzhi Li%
\thanks{%
Department of Mathematics, Southern University of Science and Technology
(SUSTech), 1088 Xueyuan Boulevard, University Town of Shenzhen, Xili,
Nanshan, Shenzhen, Guangdong Province, P.R.China (li.jz@sustech.edu.cn)} \and %
Wenlong Zhang\thanks{%
Department of Mathematics, Southern University of Science and Technology
(SUSTech), 1088 Xueyuan Boulevard, University Town of Shenzhen, Xili,
Nanshan, Shenzhen, Guangdong Province, P.R.China (zhangwl@sustech.edu.cn)} }
\date{}
\maketitle

\begin{abstract}
For the first time, we develop a convergent numerical method for the llinear
integral equation derived by M.M. Lavrent'ev in 1964 with the goal to solve
a coefficient inverse problem for a wave-like equation in 3D. The data are
non overdetermined. Convergence analysis is presented along with the
numerical results. An intriguing feature of the Lavrent'ev equation is that,
without any linearization, it reduces a highly nonlinear coefficient inverse
problem to a linear integral equation of the first kind. Nevertheless,
numerical results for that equation, which use the data generated for that
coefficient inverse problem, show a good reconstruction accuracy. This is
similar with the classical Gelfand-Levitan equation derived in 1951, which
is valid in the 1D case.
\end{abstract}



\section{Introduction}

\label{sec:1}

A member of The Russian Academy of Science Mikhail M. Lavrent'ev (1932-2010)
was one of founders of the theory of Ill-Posed and Inverse Problems. He was
a Thesis Advisor of M.V. Klibanov. In 1964, when trying to solve a 3D
Coefficient Inverse Problem (CIP) for the hyperbolic equation (\ref{2.3}),
Lavrent'ev has discovered a linear integral equation of the first kind,
which solves that CIP \cite{Lavr}, also see this equation in formula (7.18)
of the book of Lavrent'ev, Romanov and Shishatskii \cite{LRS}. The truly
intriguing point of the Lavrent'ev's idea is that even though that equation
is linear and also no linearization is made to derive it, still that
equation solves a highly nonlinear CIP for PDE (\ref{2.3}). However, the
questions about the numerical method and image quality resulting from the
numerical solution of that equation remained open. Thus, the goal of this
paper is to address the following question: \emph{Is it possible to obtain
good quality images via a numerical solution of the linear Lavrent'ev
equation using the measured data for that nonlinear CIP}? Surprisingly, the
answer is positive.

In fact, this positive answer is similar with the answer for the classical
Gelfand-Levitan (GL) linear integral equation of the second kind \cite%
{Levitan}, which was first derived in 1951 for a CIP for a 1D differential
equation \cite{GL}. Indeed, it was established numerically by Kabanikhin,
Satybev and Shishlenin in \cite{Kab1} that numerical reconstruction results
via the GL equation for the data generated by the solution of the forward
problem for that differential equation are quite accurate ones.

The above mentioned CIP is concerned with the determination of the
dielectric constant of the medium from some boundary measurements. In the
case of the backscattering data, this CIP has direct applications in the
problem of detection and identification of land mines and improvised
explosive devices, see, e.g. \cite{Khoa2,Khoa3,KlibKol} for reconstructions
from experimentally collected data. There is also an interesting problem of
the inspection of buildings using non invasive measurements of the
electromagnetic field \cite{KlibSAR,Lam}. In this case, measurements of both
backscattering and transmitted signals can be carried out, and a
corresponding CIP for equation (\ref{2.3}) is supposed to be solved.

Ramm has used the Lavrent'ev equation to prove a uniqueness theorem for a
CIP with overdetermined data \cite{Ramm}. As to the numerical studies, we
refer to the work of Bakushinskii and Leonov \cite{Bak}. In \cite{Bak} non
overdetermined data for the Lavrent'ev equation are considered. Since that
equation is a linear integral equation of the first kind, then the problem
of its solution is an ill-posed one. Hence, the numerical method of \cite%
{Bak} searches for the so-called normal solution, see, e.g. \cite{T} for
normal solutions. The computed images of \cite{Bak} are quite impressive
ones. Nevertheless, convergence of the regularized solutions to the exact
one was not proven in \cite{Bak}. We also refer to \cite{Bak1} for a similar
idea.

In this paper we consider the Lavrent'ev equation for the case of the non
overdetermined data. We construct a radically new numerical method for this
equation. For the first time, we derive from that equation a boundary value
problem (BVP) for a system of coupled elliptic PDEs with overdetermined
boundary conditions. In this derivation, we explore the orthonormal basis in
the $L_{2}$ space, which was originally introduced in \cite{KJIIP}. To solve
that BVP, we use the Quasi-Reversibility Method (QRM), which was first
introduced by Lattes and Lions \cite{LL}. Also, see \cite{KQRM} for the
recent survey on this method. In addition, we refer to \cite{Bourg1,Bourg2}
and references to these authors cited therein for some fresh ideas related
to the QRM.

To analyze our specific version of the QRM, we establish first existence and
uniqueness of the regularized solution of the QRM. Next, we turn to a much
more difficult issue of establishing the convergence rate of the regularized
solutions to the exact solution when the level of the noise in the data
tends to zero. To do this, we use the Carleman estimate of \cite{hyp}. We
recall that Carleman estimates were used in \cite{KQRM} to establish
convergence rates of regularized solutions of the QRM. However, the
technique of this paper is quite different from the one of \cite{KQRM}.
First, we do not use a cut-off function, unlike (2.15) in \cite{KQRM}.
Second, our convergence rate is in the $H^{2}$ norm rather than in the
weaker $H^{1}$ norm of \cite{KQRM}. The use of the $H^{2}$ norm is important
here since it enables us to estimate the convergence rate of resulting
solutions of the Lavrent'ev equation. Finally, our Carleman Weight Function
is significantly different from the one in section 3 of \cite{KQRM}. This
allows us to obtain convergence estimate in such a subdomain of the original
domain of interest $\Omega ,$ which is larger than a paraboloidal subdomain
of \cite{KQRM}, also see Remark 7.1 of section 7 for a relevant statement.

Since its first inception in \cite{KJIIP}, the above mentioned orthonormal
basis in $L_{2}$ was successfully used by this research group for numerical
solutions of linear inverse source problems for PDEs \cite{KN,KLN,SKN} as
well as for some versions of the globally convergent convexification method
for CIPs \cite{Khoa1,Khoa2,Khoa3,KEIT,TTTP1,TTTP2}.

We end up with presenting numerical results which solve CIPs for that
hyperbolic PDE using the Lavrent'ev equation.

In section 2 we derive the Lavrent'ev equation from a CIP\ for a hyperbolic
PDE. In section 3 we briefly review the above mentioned orthonormal basis.
In section 4, we derive that BVP\ for a system of linear coupled elliptic
PDEs. In section 5, we introduce the QRM for that BVP. In section 6, we
formulate Carleman estimates. In section 7 we provide convergence analysis
for the BVP. Numerical studies are carried out in section 8.

\section{Derivation of the Lavrent'ev Equation}

\label{sec:2}

Our derivation is different in some respects from the one of Lavrent'ev in 
\cite{Lavr}. Below $\mathbf{x}=\left( x,y,z\right) \in \mathbb{R}^{3}.$ Let
the numbers $A,b,d>0$. Let the domain $\Omega \subset \mathbb{R}^{3}$ be a
cube,%
\begin{equation}
\Omega =\left\{ \mathbf{x}=\left( x,y,z\right) :-A<x,y,z<A\right\} .
\label{2.0}
\end{equation}%
Let $L_{\text{s}}$ be a line of locations of point sources $\mathbf{x}_{0}$,
and this line is located below the cube $\Omega ,$%
\begin{equation}
L_{\text{s}}=\left\{ \mathbf{x}_{0}=\left( x_{0},0,-A-b\right) ,\text{ }%
-d<x_{0}<d\right\} .  \label{2.01}
\end{equation}%
Let the function $c\left( \mathbf{x}\right) $ be sufficiently smooth in $%
\mathbb{R}^{3}$ and also 
\begin{equation}
c\left( \mathbf{x}\right) \geq 1\text{ in }\mathbb{R}^{3},  \label{2.1}
\end{equation}%
\begin{equation}
c\left( \mathbf{x}\right) =1\text{ in }\mathbb{R}^{3}\diagdown \Omega .
\label{2.2}
\end{equation}

Consider the following Cauchy problem%
\begin{equation}
c\left( \mathbf{x}\right) u_{tt}=\Delta u,\text{ }\mathbf{x}\in \mathbb{R}%
^{3},t>0,  \label{2.3}
\end{equation}%
\begin{equation}
u\left( \mathbf{x},0\right) =0,u_{t}\left( \mathbf{x},0\right) =\delta
\left( \mathbf{x}-\mathbf{x}_{0}\right) .  \label{2.4}
\end{equation}

\textbf{Coefficient Inverse Problem (CIP).} \emph{Let }$S\subset \mathbb{R}%
^{3}$\emph{be a finite part of a plane and }$S\cap \overline{\Omega }=S\cap
L_{\text{s}}=\varnothing .$\emph{\ Find the function }$c\left( \mathbf{x}%
\right) $\emph{\ for }$\mathbf{x}\in \Omega ,$\emph{\ given the functions} $%
p_{0}\left( \mathbf{x},\mathbf{x}_{0},t\right) $ and $p_{1}\left( \mathbf{x},%
\mathbf{x}_{0},t\right) ,$ 
\begin{equation}
u\left( \mathbf{x},\mathbf{x}_{0},t\right) =p_{0}\left( \mathbf{x},\mathbf{x}%
_{0},t\right) ,\text{ }\left( \mathbf{x},\mathbf{x}_{0},t\right) \in S\times
L_{\text{s}}\times \left( 0,\infty \right) ,  \label{2.40}
\end{equation}%
\begin{equation}
\frac{\partial }{\partial n\left( \mathbf{x}\right) }u\left( \mathbf{x},%
\mathbf{x}_{0},t\right) =p_{1}\left( \mathbf{x},\mathbf{x}_{0},t\right) ,%
\text{ }\left( \mathbf{x},\mathbf{x}_{0},t\right) \in S\times L_{\text{s}%
}\times \left( 0,\infty \right) ,  \label{2.41}
\end{equation}%
\emph{where} $n=n\left( \mathbf{x}\right) $ \emph{is a unit normal vector to}
$S$ \emph{at the point} $\mathbf{x}\in S.$

We assume that the function $c\left( \mathbf{x}\right) $ satisfies the
non-trapping condition. In other words, for any bounded domain $G\subset 
\mathbb{R}^{3}$, any geodesic line generated by $c\left( \mathbf{x}\right) $
and starting in $G$ leaves this domain in a finite time, see the book of
Vainberg \cite{V2}. In this case the function $u\left( \mathbf{x},\mathbf{x}%
_{0},t\right) $ decays exponentially as $t\rightarrow \infty $ together with
its derivatives with respect to $\mathbf{x,}t$ up to the second order, and
this decay is uniform with respect to $\mathbf{x}\in G$ for any bounded
domain $G\subset \mathbb{R}^{3},$ see Lemma 6 of Chapter 10 of the book \cite%
{V2} and Remark 3 after this lemma.

Relaying on these results, apply to the function $u$ the Fourier transform
with respect to $t$, 
\begin{equation}
w\left( \mathbf{x},\mathbf{x}_{0},k\right) =\mathop{\displaystyle \int}%
\limits_{0}^{\infty }u\left( \mathbf{x},\mathbf{x}_{0},t\right)
e^{ikt}dt,k\geq 0.  \label{2.6}
\end{equation}%
Due to that exponential decay, one can differentiate the function $w\left( 
\mathbf{x},\mathbf{x}_{0},k\right) $ at least twice with respect to $\mathbf{%
x}$ and infinitely many times with respect to $k$. Furthermore, these
derivatives can be calculated under the integral sign. A combination of
Theorem 3.3 of \cite{V1} and Theorem 6 of Chapter 9 of \cite{V2} ensures
that 
\begin{equation}
\Delta _{\mathbf{x}}w+k^{2}w+k^{2}\left( c\left( \mathbf{x}\right) -1\right)
w=-\delta \left( \mathbf{x}-\mathbf{x}_{0}\right) ,  \label{2.7}
\end{equation}%
\begin{equation}
\frac{\partial w}{\partial r}-ikw=o\left( \frac{1}{r}\right) ,\text{ }%
r\rightarrow \infty ,r=\left\vert \mathbf{x}-\mathbf{x}_{0}\right\vert .
\label{2.8}
\end{equation}%
Denote 
\begin{equation}
w_{0}\left( \mathbf{x}-\mathbf{x}_{0},k\right) =\frac{\exp \left(
ik\left\vert \mathbf{x}-\mathbf{x}_{0}\right\vert \right) }{4\pi \left\vert 
\mathbf{x}-\mathbf{x}_{0}\right\vert }.  \label{2.80}
\end{equation}%
It is well known that the function $w\left( \mathbf{x},\mathbf{x}%
_{0},k\right) $ satisfies the Lippmann-Schwinger integral equation \cite%
{Colton}%
\begin{equation}
w\left( \mathbf{x},\mathbf{x}_{0},k\right) =w_{0}\left( \mathbf{x}-\mathbf{x}%
_{0},k\right) +k^{2}\mathop{\displaystyle \int}\limits_{\Omega }w_{0}\left( 
\mathbf{x}-\mathbf{\xi },k\right) \left( c\left( \mathbf{\xi }\right)
-1\right) w\left( \mathbf{\xi },x_{0},k\right) d\mathbf{\xi .}  \label{2.9}
\end{equation}%
Differentiate (\ref{2.9}) twice with respect to $k$ and then set $k=0.$
Denote%
\begin{equation}
q\left( \mathbf{x}\right) =c\left( \mathbf{x}\right) -1.  \label{2.10}
\end{equation}%
\begin{equation}
w^{0}\left( \mathbf{x},\mathbf{x}_{0}\right) =\partial _{k}^{2}w\left( 
\mathbf{x},\mathbf{x}_{0},k\right) \mid _{k=0}.  \label{2.11}
\end{equation}%
By (\ref{2.9}) and (\ref{2.10}) 
\begin{equation}
w^{0}\left( \mathbf{x},\mathbf{x}_{0}\right) =-\mathop{\displaystyle \int}%
\limits_{0}^{\infty }t^{2}u\left( \mathbf{x},\mathbf{x}_{0},t\right) dt.
\label{2.12}
\end{equation}

It follows from (\ref{2.80})-(\ref{2.11}) that%
\begin{equation}
\frac{1}{8\pi ^{2}}\mathop{\displaystyle \int}\limits_{\Omega }\frac{q\left( 
\mathbf{\xi }\right) }{\left\vert \mathbf{x}-\mathbf{\xi }\right\vert
\left\vert \xi -\mathbf{x}_{0}\right\vert }d\mathbf{\xi =}w^{0}\left( 
\mathbf{x},\mathbf{x}_{0}\right) -\frac{\left\vert \mathbf{x}-\mathbf{x}%
_{0}\right\vert }{4\pi }.  \label{2.13}
\end{equation}%
Equation (\ref{2.13}) with respect to the function $q\left( \mathbf{x}%
\right) $ is the Lavrent'ev equation.

We now reformulate the CIP (\ref{2.40}), (\ref{2.41}) in terms of equation (%
\ref{2.13}). Denote 
\begin{equation}
f_{0}\left( \mathbf{x},\mathbf{x}_{0}\right) =-\mathop{\displaystyle \int}%
\limits_{0}^{\infty }p_{0}\left( \mathbf{x},\mathbf{x}_{0},t\right) t^{2}dt-%
\frac{\left\vert \mathbf{x}-\mathbf{x}_{0}\right\vert }{4\pi },\text{ for }%
\left( \mathbf{x},\mathbf{x}_{0}\right) \in S\times L_{\text{s}},
\label{2.14}
\end{equation}%
\begin{equation}
f_{1}\left( \mathbf{x},\mathbf{x}_{0}\right) =-\mathop{\displaystyle \int}%
\limits_{0}^{\infty }p_{1}\left( \mathbf{x},\mathbf{x}_{0},t\right) t^{2}dt-%
\frac{\partial }{\partial n\left( \mathbf{x}\right) }\left( \frac{\left\vert 
\mathbf{x}-\mathbf{x}_{0}\right\vert }{4\pi }\right) ,\text{ for }\left( 
\mathbf{x},\mathbf{x}_{0}\right) \in S\times L_{\text{s}}.  \label{2.15}
\end{equation}%
Denote 
\begin{equation}
g_{0}\left( \mathbf{x},\mathbf{x}_{0}\right) =2\pi f_{0}\left( \mathbf{x},%
\mathbf{x}_{0}\right) -\frac{\left\vert \mathbf{x}-\mathbf{x}_{0}\right\vert 
}{2},\text{ for }\left( \mathbf{x},\mathbf{x}_{0}\right) \in S\times L_{%
\text{s}},  \label{2.16}
\end{equation}%
\begin{equation}
g_{1}\left( \mathbf{x},\mathbf{x}_{0}\right) =2\pi f_{1}\left( \mathbf{x},%
\mathbf{x}_{0}\right) -\frac{\partial }{\partial n\left( \mathbf{x}\right) }%
\left( \frac{\left\vert \mathbf{x}-\mathbf{x}_{0}\right\vert }{2}\right) ,%
\text{ for }\left( \mathbf{x},\mathbf{x}_{0}\right) \in S\times L_{\text{s}}.
\label{2.17}
\end{equation}%
Then (\ref{2.13}) implies: 
\begin{equation}
\frac{1}{4\pi }\mathop{\displaystyle \int}\limits_{\Omega }\frac{q\left( 
\mathbf{\xi }\right) }{\left\vert \mathbf{x}-\mathbf{\xi }\right\vert
\left\vert \xi -\mathbf{x}_{0}\right\vert }d\mathbf{\xi =}g_{0}\left( 
\mathbf{x},\mathbf{x}_{0}\right) ,\text{ for }\left( \mathbf{x},\mathbf{x}%
_{0}\right) \in S\times L_{\text{s}},  \label{2.19}
\end{equation}%
\begin{equation}
\frac{\partial }{\partial n\left( \mathbf{x}\right) }\left( \frac{1}{4\pi }%
\mathop{\displaystyle \int}\limits_{\Omega }\frac{q\left( \mathbf{\xi }%
\right) }{\left\vert \mathbf{x}-\mathbf{\xi }\right\vert \left\vert \xi -%
\mathbf{x}_{0}\right\vert }d\mathbf{\xi }\right) \mathbf{=}g_{1}\left( 
\mathbf{x},\mathbf{x}_{0}\right) ,\text{ for }\left( \mathbf{x},\mathbf{x}%
_{0}\right) \in S\times L_{\text{s}}.  \label{2.20}
\end{equation}

Therefore, we have reduced the CIP (\ref{2.40}), (\ref{2.41}) to the
following inverse problem:

\textbf{Problem For The Lavrent'ev Equation to Solve.} Find the function $%
q\left( \mathbf{x}\right) ,\mathbf{x}\in \Omega $ from functions $%
g_{0}\left( \mathbf{x},\mathbf{x}_{0}\right) $ and $g_{1}\left( \mathbf{x},%
\mathbf{x}_{0}\right) $ given in (\ref{2.19}), (\ref{2.20}).

\textbf{Remarks 2.1}:

1. \emph{Note that this problem is non overdetermined. Indeed, in this
problem }$m=n=3$\emph{, where }$m$\emph{\ is the number of free variables in
the data and }$n$\emph{\ is the number of free variables in the unknown
function} $q\left( \mathbf{x}\right) .$ \emph{On the other hand, }$m>n$\emph{%
\ in the overdetermined case.}

\emph{2. It follows from (\ref{2.19}) that the right hand side of this
equation is a harmonic function outside of the domain }$\Omega $\emph{. If,
for example }$S=\partial \Omega ,$\emph{\ then this function can be uniquely
determined outside }$\Omega $\emph{\ by its boundary values }$g_{0}\left( 
\mathbf{x},\mathbf{x}_{0}\right) $\emph{\ at }$S.$\emph{\ Hence, in this
case we do not need to assign the normal derivative }$g_{1}\left( \mathbf{x},%
\mathbf{x}_{0}\right) $\emph{\ at }$S$\emph{\ since it will be calculated
using }$g_{0}\left( \mathbf{x},\mathbf{x}_{0}\right) .$\emph{\ However,
since }$S\neq \partial \Omega $\emph{\ above, then we indeed need to assign }%
$g_{1}\left( \mathbf{x},\mathbf{x}_{0}\right) $\emph{\ in (\ref{2.20}).}

\emph{3. Even though we have sometimes worked above with complex valued
functions, since we will work from now on with equations (\ref{2.19}), (\ref%
{2.20}) with a real valued function }$q$\emph{, then all functions below are
real valued ones. Also, all function spaces below are spaces of real valued
functions.}

\section{A Special Orthonormal Basis}

\label{sec:3}

Recall that by (\ref{2.01}) the point source $\mathbf{x}_{0}=\left(
x_{0},0,-A-b\right) \in L_{\text{s}}$ and $x_{0}\in \left( -d,d\right) .$
Prior to the description of our numerical method, we present in this section
a special orthonormal basis $\left\{ \Psi _{n}\left( x_{0}\right) \right\}
_{n=0}^{\infty }$ in $L_{2}\left( -d,d\right) .$ This basis has the
following two properties:

\begin{enumerate}
\item $\Psi _{n}\in C^{1}\left[ -d,d\right] ,$ $\forall n=1,2,\dots $

\item Let $\left\{ ,\right\} $ be the scalar product in $L_{2}\left(
-d,d\right) $. Denote $a_{mn}=\left\{ \Psi _{n}^{\prime },\Psi _{m}\right\}
. $ Then the matrix $M_{N}=\left( a_{mn}\right) _{m,n=0}^{N}$ is invertible
for any $N=1,2,\dots $
\end{enumerate}

Note that if one would use either the basis of standard orthonormal
polynomials or the trigonometric basis, then the first derivative of the
first element of either of them would be identically zero. Hence, the matrix 
$M_{N}$ would not be invertible in this case.

For the convenience of the reader, we now describe that basis. Consider the
set of functions $\left\{ \xi _{n}\left( x_{0}\right) \right\}
_{n=0}^{\infty }=\{x_{0}^{n}e^{x_{0}}\}_{n=0}^{\infty }.$ This is a set of
linearly independent functions. Besides, this set is complete in the space $%
L_{2}(-d,d)$. After applying the Gram-Schmidt orthonormalization procedure
to this set, we obtain the orthonormal basis $\{\Psi _{n}\left( x_{0}\right)
\}_{n=1}^{\infty }$ in $L_{2}(-d,d)$. In fact, the function $\Psi
_{n}(x_{0}) $ has the form $\Psi _{n}(x_{0})=P_{n}(x_{0})e^{x_{0}}$, $%
\forall n\geq 1,$ where $P_{n}$ is a polynomial of the degree $n$. Thus,
functions $\Psi _{n}(x_{0})$ are polynomials which are orthonormal in $%
L_{2}(-d,d)$ with the weight $e^{2x_{0}}.$ The matrix $M_{N}$ is invertible
since $\det M_{N}=1.$ Indeed elements $a_{mn}=\left( \Psi _{n}^{\prime
},\Psi _{m}\right) $ of this matrix are such that \cite{KJIIP}%
\[
a_{mn}=\left\{ 
\begin{array}{c}
1\text{ if }m=n\text{,} \\ 
0\text{ if }m>n.%
\end{array}%
\right. 
\]%
In principle, the Gram-Schmidt orthonormalization procedure is unstable.
However, in all our above referenced past works where the above basis was
used we have always managed to find such a number $N$ that it was stable for
this $N$. This means that $N$ is the regularization parameter of this
procedure.

\textbf{Remark 3.1}. \emph{As soon as }$N$\emph{\ is fixed, we work within
the framework of an approximate mathematical model. Estimating convergence
of our method at }$N\rightarrow \infty $\emph{\ is a very challenging
problem which we cannot address at the moment}$.$\emph{The underlying reason
of this is the ill-posed nature of problem (\ref{2.19}), (\ref{2.20}).
Still, our numerical reconstructions of section 8 are of a good quality. We
point out that the similar approximate mathematical models are quite often
used in the field of inverse and ill-posed problems without proofs of
convergence at  }$N\rightarrow \infty ,$ \emph{again because of the
ill-posed nature of corresponding problems. In this regard we refer to, e.g.
works \cite{G,Kab1,Kab2,Kab3} of the authors, who are not members of this
research team. Numerical results in all these publications are quite good
ones. More detailed discussions of this issue can be found in \cite%
{Khoa1,hyp}. Even though similar truncations were made in these \cite%
{Khoa1,hyp}, still good quality images were obtained.}

\section{Numerical Method for Problem (\protect\ref{2.19}), (\protect\ref%
{2.20})}

\label{sec:4}

Below for each Banach space $Z$ and for each integer $m$, we denote $Z_{m}$
the space of $m-$dimensional vectors $D=\left( D_{1},...,D_{m}\right)
,D_{j}\in Z$ with the norm%
\[
\left\Vert D\right\Vert _{Z_{m}}=\left( \mathop{\displaystyle \sum }%
\limits_{j=1}^{m}\left\Vert D_{j}\right\Vert _{Z}^{2}\right) ^{1/2}. 
\]%
Keeping in mind locations of point sources in (\ref{2.01}), we assume below
for the surface $S$ that it is either the backscattering part of the
boundary of the cube $\Omega $ in (\ref{2.0})$,S=S_{bsc}$ or its transmitted
part $S=S_{tr}$, where 
\begin{equation}
S_{bsc}=\left\{ \mathbf{x}=\left( x,y,z\right) :-A<x,y<A,z=-A\right\} ,
\label{4.1}
\end{equation}%
\begin{equation}
S_{tr}=\left\{ \mathbf{x}=\left( x,y,z\right) :-A<x,y<A,z=A\right\} .
\label{4.2}
\end{equation}%
Denote%
\[
v\left( \mathbf{x},\mathbf{x}_{0}\right) =2\pi \left[ w^{0}\left( \mathbf{x},%
\mathbf{x}_{0}\right) -\frac{\left\vert \mathbf{x}-\mathbf{x}_{0}\right\vert 
}{4\pi }\right] . 
\]%
By (\ref{2.13})%
\begin{equation}
v\left( \mathbf{x},\mathbf{x}_{0}\right) =\frac{1}{4\pi }\mathop{%
\displaystyle \int}\limits_{\Omega }\frac{q\left( \mathbf{\xi }\right) }{%
\left\vert \mathbf{x}-\mathbf{\xi }\right\vert \left\vert \xi -\mathbf{x}%
_{0}\right\vert }d\mathbf{\xi }.  \label{3.1}
\end{equation}%
Hence,%
\begin{equation}
\Delta _{\mathbf{x}}v\left( \mathbf{x},\mathbf{x}_{0}\right) =-\frac{q\left( 
\mathbf{x}\right) }{\left\vert \mathbf{x}-\mathbf{x}_{0}\right\vert },\text{ 
}\mathbf{x}\in \Omega ,\mathbf{x}_{0}\in L_{\text{s}}.  \label{3.2}
\end{equation}%
Also, by (\ref{2.19}) and (\ref{2.20})%
\begin{equation}
v\left( \mathbf{x},\mathbf{x}_{0}\right) =g_{0}\left( \mathbf{x},\mathbf{x}%
_{0}\right) ,\text{ }\mathbf{x}\in S,\mathbf{x}_{0}\in L_{\text{s}},
\label{3.02}
\end{equation}%
\begin{equation}
\frac{\partial }{\partial n\left( \mathbf{x}\right) }v\left( \mathbf{x},%
\mathbf{x}_{0}\right) =g_{1}\left( \mathbf{x},\mathbf{x}_{0}\right) ,\text{ }%
\mathbf{x}\in S,\mathbf{x}_{0}\in L_{\text{s}}.  \label{3.03}
\end{equation}%
Multiply both sides of (\ref{3.2}) by $\left\vert \mathbf{x}-\mathbf{x}%
_{0}\right\vert $ and, slightly abusing one of notations of section 2,
denote 
\begin{equation}
u\left( \mathbf{x},\mathbf{x}_{0}\right) =v\left( \mathbf{x},\mathbf{x}%
_{0}\right) \left\vert \mathbf{x}-\mathbf{x}_{0}\right\vert ,\text{ }\mathbf{%
x}\in \Omega ,\mathbf{x}_{0}\in L_{\text{s}}.  \label{3.3}
\end{equation}%
By (\ref{3.2}) 
\begin{equation}
\left\vert \mathbf{x}-\mathbf{x}_{0}\right\vert \Delta _{\mathbf{x}}v\left( 
\mathbf{x},\mathbf{x}_{0}\right) =-q\left( \mathbf{x}\right) ,\text{ }%
\mathbf{x}\in \Omega ,\mathbf{x}_{0}\in L_{\text{s}}.  \label{3.4}
\end{equation}%
By (\ref{3.3}) the left hand side of (\ref{3.4}) can be rewritten as%
\[
\left\vert \mathbf{x}-\mathbf{x}_{0}\right\vert \Delta _{\mathbf{x}}v=\Delta
_{\mathbf{x}}\left( \left\vert \mathbf{x}-\mathbf{x}_{0}\right\vert v\right)
-2\nabla \left( \left\vert \mathbf{x}-\mathbf{x}_{0}\right\vert \right)
\nabla v-\Delta \left( \left\vert \mathbf{x}-\mathbf{x}_{0}\right\vert
\right) v 
\]%
\begin{equation}
=\Delta _{\mathbf{x}}u-2\frac{\mathbf{x}-\mathbf{x}_{0}}{\left\vert \mathbf{x%
}-\mathbf{x}_{0}\right\vert ^{2}}\nabla _{\mathbf{x}}u+\frac{4}{\left\vert 
\mathbf{x}-\mathbf{x}_{0}\right\vert ^{2}}u.  \label{3.40}
\end{equation}%
To improve the stability of our method, we set:%
\begin{equation}
\frac{\partial }{\partial n\left( \mathbf{x}\right) }u\left( \mathbf{x},%
\mathbf{x}_{0}\right) =\left\{ 
\begin{array}{c}
0\text{ for }\mathbf{x}\in S_{tr}\text{ if }S=S_{bsc}, \\ 
0\text{ for }\mathbf{x}\in S_{bsc}\text{ if }S=S_{tr}, \\ 
0\text{ for }\mathbf{x}\in \left\{ x,y=\pm A,z\in \left( -A,A\right)
\right\} ,%
\end{array}%
\right. \mathbf{x}_{0}\in L_{\text{s}}.  \label{3.50}
\end{equation}%
It follows from (\ref{4.1}), (\ref{4.2}) and (\ref{3.1}) that conditions (%
\ref{3.50}) are approximately satisfied if the number $A$ is sufficiently
large and the support of $q\left( \mathbf{x}\right) $ is far from the
boundary $\partial \Omega .$\ 

Next, (\ref{2.01}), (\ref{3.3}) and (\ref{3.40}) imply that (\ref{3.2}) is
equivalent with%
\begin{equation}
\Delta _{\mathbf{x}}u-\frac{2}{\left( x-x_{0}\right) ^{2}+y^{2}+z^{2}}\left(
\left( x-x_{0}\right) u_{x}+yu_{y}+zu_{z}-2u\right) =-\frac{q\left(
x,y,z\right) }{2\pi },  \label{3.6}
\end{equation}%
where $\left( x,y,z\right) \in \Omega ,x_{0}\in \left( -d,d\right) .$
Differentiate both sides of (\ref{3.6}) with respect to $x_{0}.$ We obtain
for $\left( x,y,z\right) \in \Omega ,x_{0}\in \left( -d,d\right) :$%
\begin{equation}
\Delta _{\mathbf{x}}u_{x_{0}}-2\frac{\partial }{\partial x_{0}}\left[ \frac{%
\left( x-x_{0}\right) u_{x}+yu_{y}+zu_{z}-2u}{\left( x-x_{0}\right)
^{2}+y^{2}+z^{2}}\right] =0.  \label{3.7}
\end{equation}

We now represent the function $u\left( \mathbf{x},\mathbf{x}_{0}\right)
=u\left( x,y,z,x_{0}\right) $ as a truncated Fourier series with respect to
the basis $\left\{ \Psi _{n}\left( x_{0}\right) \right\} ,$%
\begin{equation}
u\left( x,y,z,x_{0}\right) =\mathop{\displaystyle \sum }%
\limits_{n=1}^{N}u_{n}\left( x,y,z\right) \Psi _{n}\left( x_{0}\right) ,
\label{3.8}
\end{equation}%
\begin{equation}
u_{x_{0}}\left( x,y,z,x_{0}\right) =\mathop{\displaystyle \sum }%
\limits_{n=1}^{N}u_{n}\left( x,y,z\right) \Psi _{n}^{\prime }\left(
x_{0}\right) ,  \label{3.9}
\end{equation}%
for $\left( x,y,z\right) \in \Omega ,x_{0}\in \left( -d,d\right) .$ Here the
number $N>1$ should be chosen numerically. Therefore, coefficients $%
u_{n}\left( x,y,z\right) $ of the Fourier expansion (\ref{3.8}) are unknown
now, and we should focus on the search for these functions.

Substituting (\ref{3.8}) and (\ref{3.9}) in (\ref{3.7}), we obtain for $%
\mathbf{x}\in \Omega ,x_{0}\in \left( -d,d\right) $ 
\begin{equation}
\mathop{\displaystyle \sum }\limits_{n=1}^{N}\Delta _{\mathbf{x}}u_{n}\left( 
\mathbf{x}\right) \Psi _{n}^{\prime }\left( x_{0}\right)  \label{3.10}
\end{equation}%
\[
-\frac{\partial }{\partial x_{0}}\left\{ \frac{2}{\left\vert \mathbf{x}-%
\mathbf{x}_{0}\right\vert ^{2}}\mathop{\displaystyle \sum }\limits_{n=1}^{N}%
\left[ \left( \mathbf{x}-\mathbf{x}_{0}\right) \nabla u_{n}\left( \mathbf{x}%
\right) -2u_{n}\left( \mathbf{x}\right) \right] \Psi _{n}\left( x_{0}\right)
\right\} =0. 
\]%
Multiply both sides of (\ref{3.10}) sequentially by functions $\Psi
_{m}\left( x_{0}\right) ,m=1,...,N$ and integrate then with respect to $%
x_{0}\in \left( -d,d\right) .$ Denote the $n-$D vector of unknown
coefficients in (\ref{3.8}) as $U\left( \mathbf{x}\right) =\left(
u_{1},...u_{N}\right) ^{T}\left( \mathbf{x}\right) .$ We obtain a system of
coupled elliptic equations of the second order,%
\begin{equation}
L\left( U\right) =\Delta U+A_{1}\left( \mathbf{x}\right) U_{x}+A_{2}\left( 
\mathbf{x}\right) U_{y}+A_{3}\left( \mathbf{x}\right) U_{z}+A_{0}\left( 
\mathbf{x}\right) U=0,  \label{3.11}
\end{equation}%
where $A_{j}\left( \mathbf{x}\right) \in C_{N^{2}}\left( \overline{\Omega }%
\right) $ are $N\times N$ matrices. It follows from (\ref{3.02}) and (\ref%
{3.03}) that the following boundary vector functions $\varphi _{0}\left( 
\mathbf{x}\right) ,\varphi _{1}\left( \mathbf{x}\right) $ are known:%
\begin{equation}
U\left( \mathbf{x}\right) =\varphi _{0}\left( \mathbf{x}\right) ,\text{ }%
\frac{\partial }{\partial n\left( \mathbf{x}\right) }U\left( \mathbf{x}%
\right) =\varphi _{1}\left( \mathbf{x}\right) \text{ for }\mathbf{x}\in S.
\label{3.12}
\end{equation}%
In addition (\ref{3.50}) implies that%
\begin{equation}
\frac{\partial }{\partial n\left( \mathbf{x}\right) }U\left( \mathbf{x}%
\right) =0\text{ for }\mathbf{x}\in \left\{ \left\vert x\right\vert
,\left\vert y\right\vert =A,z\in \left( -A,A\right) \right\} ,  \label{3.13}
\end{equation}%
\begin{equation}
\frac{\partial }{\partial n\left( \mathbf{x}\right) }U\left( \mathbf{x}%
\right) =0\text{ for }\mathbf{x}\in S_{tr}\text{ if }S=S_{bsc},  \label{3.14}
\end{equation}%
\begin{equation}
\frac{\partial }{\partial n\left( \mathbf{x}\right) }U\left( \mathbf{x}%
\right) =0\text{ for }\mathbf{x}\in S_{bsc}\text{ if }S=S_{tr}\text{ }.
\label{3.15}
\end{equation}%
Thus, depending on the location of the surface $S$: whether $S=S_{bsc}$ or $%
S=S_{tr},$ we actually got two boundary value problems:

\textbf{Boundary Value Problem 1 (BVP1)}. \emph{Solve system (\ref{3.11})
with the boundary conditions (\ref{3.12})-(\ref{3.14}).}

\textbf{Boundary Value Problem 2 (BVP2)}. \emph{Solve system (\ref{3.11})
with the boundary conditions (\ref{3.12}), (\ref{3.13}), (\ref{3.15}).}

\textbf{Remark 4.1.} \emph{Since we have the Cauchy data at the surface }$S$%
\emph{\ in (\ref{3.12}) for the elliptic system (\ref{3.11}), then each of
these BVPs has at most one solution. This follows from the classical unique
continuation theorem for elliptic PDEs, see, e.g. Chapter \ 4 of \cite{LRS}. 
}

\section{The Quasi-Reversibility Method For Problems 1,2}

\label{5}

Both these problems are overdetermined ones due to the overdetermination of
the boundary data (\ref{3.12}). Thus, we solve both problems via the
Quasi-Reversibility Method (QRM), which is well suitable for solving
overdetermined BVPs for PDEs. Following the survey \cite{KQRM} on the QRM,
we consider two minimization problems. We remind that a significant
difference with \cite{KQRM} is in the convergence Theorem 7.2.

Denote $\partial _{1}\Omega =\partial \Omega \diagdown \left\{ z=A\right\}
,\partial _{2}\Omega =\partial \Omega \diagdown \left\{ z=-A\right\} .$
Introduce two subspaces of the space $H^{3}\left( \Omega \right) ,$%
\begin{equation}
H_{0,bsc}^{3}\left( \Omega \right) =\left\{ u\in H^{3}\left( \Omega \right)
:u\mid _{S_{bsc}}=\partial _{n\left( \mathbf{x}\right) }u\mid _{\partial
_{1}\Omega }=0\right\} ,  \label{5.1}
\end{equation}%
\[
H_{0,tr}^{3}\left( \Omega \right) =\left\{ u\in H^{3}\left( \Omega \right)
:u\mid _{S_{tr}}=\partial _{n\left( \mathbf{x}\right) }u\mid _{\partial
_{2}\Omega }=0\right\} . 
\]

To solve BVPs 1 and 2, we consider Minimization Problems 1 and 2
respectively. Below $\mathbf{x}=\left( x,y,z\right) .$ We assume that, in
the case of BVP 1, there exists such a vector function $F_{1}\in
H_{N}^{3}\left( \Omega \right) $ that satisfies boundary conditions (\ref%
{3.12})-(\ref{3.14}). And in the case of Problem 2, we assume the existence
of such a function $F_{2}\in H_{N}^{3}\left( \Omega \right) $ that satisfies
boundary conditions (\ref{3.12}), (\ref{3.13}), (\ref{3.15}). Denote%
\begin{equation}
\left\{ 
\begin{array}{c}
V=U-F_{1}\text{ in the case of BVP1,} \\ 
W=U-F_{2}\text{ in the case of BVP2.}%
\end{array}%
\right.  \label{5.2}
\end{equation}%
Hence, $V\in H_{0,bsc,N}^{3}\left( \Omega \right) $ and $W\in
H_{0,tr,N}^{3}\left( \Omega \right) .$

\textbf{Minimization Problem 1.} Let $L$ be the operator in (\ref{3.11}).
Minimize the functional $J_{\gamma }$ on the space $H_{0,bsc,N}^{3}\left(
\Omega \right) ,$%
\begin{equation}
J_{\gamma }\left( V\right) =\mathop{\displaystyle \int}\limits_{\Omega
}\left\vert L\left( V+F_{1}\right) \right\vert ^{2}d\mathbf{x+}\gamma
\left\Vert V\right\Vert _{H_{N}^{3}\left( \Omega \right) }^{2}.
\label{5.200}
\end{equation}%
Here $\gamma \in \left( 0,1\right) $ is the regularization parameter.

\textbf{Minimization Problem 2}. Minimize the functional $I_{\gamma }$ on
the space $H_{0,tr,N}^{3}\left( \Omega \right) $ 
\begin{equation}
I_{\gamma }\left( W\right) =\mathop{\displaystyle \int}\limits_{\Omega
}\left\vert L\left( W+F_{2}\right) \right\vert ^{2}d\mathbf{x+}\gamma
\left\Vert W\right\Vert _{H_{N}^{3}\left( \Omega \right) }^{2}.
\label{5.201}
\end{equation}

\textbf{Remark 5.1}. \emph{As soon as a minimizer of either of these
functionals is found, we reconstruct the vector function }$U$\emph{\ by
either of formulas (\ref{5.2}). Then we reconstruct the function }$u\left(
x,y,z,x_{0}\right) $\emph{\ via (\ref{3.8}) and then reconstruct the target
function }$q\left( x,y,z\right) $\emph{\ via (\ref{3.6}). In our
computations we use }$x_{0}=0$\emph{\ to reconstruct }$q$\emph{\ via (\ref%
{3.6}).}

\section{Carleman Estimates}

\label{6}

Let $\lambda \geq 1$ be a parameter, which is always used in \ CWFs and let
the number $R>A$. Introduce two CWFs $\mu _{1,\lambda }\left( z\right) $ and 
$\mu _{2,\lambda }\left( z\right) $ for the Minimization Problems 1 and 2
respectively,%
\begin{equation}
\mu _{1,\lambda }\left( z\right) =\exp \left( \lambda \left( z-R\right)
^{2}\right) ,z\in \left( -A,A\right) ,  \label{6.3}
\end{equation}%
\begin{equation}
\mu _{2,\lambda }\left( z\right) =\exp \left( \lambda \left( z+R\right)
^{2}\right) ,z\in \left( -A,A\right) .  \label{6.4}
\end{equation}%
Hence, $\mu _{1,\lambda }^{\prime }\left( z\right) <0$ for $z\in \left[ -A,A%
\right] $ and $\mu _{2,\lambda }^{\prime }\left( z\right) >0$ for $z\in %
\left[ -A,A\right] .$ This means that the function $\mu _{1,\lambda }\left(
z\right) $ is decreasing for $z\in \left[ -A,A\right] $ and the function $%
\mu _{2,\lambda }\left( z\right) $ is increasing for $z\in \left[ -A,A\right]
.$ Below $C=C\left( A,R\right) >0$ denotes different constants depending
only on the number $A$ in (\ref{2.0}). Also, $B=B\left( A,R,L\right) >0$
denotes different constants depending only on numbers $A,R$ and the maximum
of $C_{N^{2}}\left( \overline{\Omega }\right) -$norms of matrices of the
operator $L$ in (\ref{3.11}).

\textbf{Theorem 6.1} (Carleman estimate)\emph{\ There exists constants }$%
\lambda _{0}=\lambda _{0}\left( A\right) >1$\emph{\ and }$C=C\left(
A,R\right) >0$\emph{\ such that the following Carleman estimate holds for
all functions }$u\in H_{0,bsc}^{3}\left( \Omega \right) $\emph{\ and for all 
}$\lambda \geq \lambda _{0}:$%
\begin{equation}
\mathop{\displaystyle \int}\limits_{\Omega }\left( \Delta u\right) ^{2}\mu
_{1,\lambda }^{2}d\mathbf{x\geq }\frac{C}{\lambda }\mathop{\displaystyle
\sum }\limits_{i,j=1}^{3}\mathop{\displaystyle \int}\limits_{\Omega
}u_{x_{i}x_{j}}^{2}\mu _{1,\lambda }^{2}d\mathbf{x+}C\lambda %
\mathop{\displaystyle \int}\limits_{\Omega }\left( \left\vert \nabla
u\right\vert ^{2}+\lambda ^{2}u^{2}\right) \mu _{1,\lambda }^{2}d\mathbf{x}
\label{6.5}
\end{equation}%
\[
-C\lambda ^{3}e^{2\lambda \left( A-R\right) ^{2}}\left\Vert u\right\Vert
_{H^{3}\left( \Omega \right) }^{2}. 
\]%
\textbf{Theorem 6.2} (Carleman estimate). \emph{There exists constants }$%
\lambda _{0}=\lambda _{0}\left( A\right) >0$\emph{\ and }$C=C\left(
A,R\right) >0$\emph{\ such that the following Carleman estimate holds for
all functions }$u\in H_{0,tr}^{3}\left( \Omega \right) $\emph{\ and for all }%
$\lambda \geq \lambda _{0}:$%
\[
\mathop{\displaystyle \int}\limits_{\Omega }\left( \Delta u\right) ^{2}\mu
_{2,\lambda }^{2}d\mathbf{x\geq }\frac{C}{\lambda }%
\mathop{\displaystyle
\sum }\limits_{i,j=1}^{3}\mathop{\displaystyle \int}\limits_{\Omega
}u_{x_{i}x_{j}}^{2}\mu _{2,\lambda }^{2}d\mathbf{x+}C\lambda %
\mathop{\displaystyle \int}\limits_{\Omega }\left( \left\vert \nabla
u\right\vert ^{2}+\lambda ^{2}u^{2}\right) \mu _{2,\lambda }^{2}d\mathbf{x} 
\]%
\[
-C\lambda ^{3}e^{2\lambda \left( A+R\right) ^{2}}\left\Vert u\right\Vert
_{H^{3}\left( \Omega \right) }^{2}. 
\]

The proof of each of these theorems is a slight modification of the proof of
Theorem 4.1 of our work \cite{hyp}, and then this modification should be
combined with the trace theorem. Hence, we do not prove Theorems 6.1 and 6.2
here.

\section{Convergence Analysis of the QRM}

\label{7}

For brevity, we provide in this section convergence analysis only for the
Minimization Problem 1, since this analysis is completely similar for the
Minimization Problem 2. The only difference is that we need to use Theorem
6.2 instead of Theorem 6.1 in the second case.

\subsection{Existence and uniqueness of the minimizer}

\label{7.1}

\textbf{Theorem 7.1.} \emph{Consider the Minimization Problem 1. For any }$%
\gamma >0$\emph{\ there exists unique minimizer }$V_{\min }\in
H_{0,bsc,N}^{3}\left( \Omega \right) $\emph{\ of the functional }$J_{\gamma
}\left( V\right) $\emph{\ and the following estimate holds}%
\[
\left\Vert V_{\min }\right\Vert _{H_{N}^{3}\left( \Omega \right) }\leq \frac{%
C}{\sqrt{\gamma }}\left\Vert F_{1}\right\Vert _{H_{N}^{3}\left( \Omega
\right) }. 
\]

The minimizer of Theorem 7.1 is called \textquotedblleft regularized
solution" in the regularization theory \cite{BK,T}. This theorem follows
immediately from Theorem 2.4 of \cite{KQRM}. That theorem was proved using
the variational principle and the Riesz theorem. However, to establish
convergence rate of the regularized solutions to the exact one when the
level of noise in the data tends to zero, we need to apply in the next
subsection the Carleman estimate of Theorem 6.1.

\subsection{Convergence rate of regularized solutions}

\label{7.2}

Denote $\left( ,\right) $ and $\left[ ,\right] $ scalar products in spaces $%
L_{2}\left( \Omega \right) $ and $H^{3}\left( \Omega \right) $ respectively.
As it is often the case in the theory of ill-posed problems, we assume that
there exists exact solution $U^{\ast }\in H_{N}^{3}\left( \Omega \right) $
of equation (\ref{3.11}) with noiseless boundary data $\varphi _{0}^{\ast
}\left( \mathbf{x}\right) ,\varphi _{1}^{\ast }\left( \mathbf{x}\right) $ in
(\ref{3.12}) and satisfying additional boundary conditions (\ref{3.13})-(\ref%
{3.15}). Since such a vector function $U^{\ast }$ exists, then there also
exists a vector function $F_{1}^{\ast }\in H_{N}^{3}\left( \Omega \right) $
satisfying these noiseless boundary conditions (\ref{3.12}) as well as
additional zero Neumann boundary conditions (\ref{3.13})-(\ref{3.15}).

Denote 
\begin{equation}
V^{\ast }=U^{\ast }-F_{1}^{\ast }.  \label{7.300}
\end{equation}
Then $V^{\ast }\in H_{0,bsc,N}^{3}\left( \Omega \right) $ and 
\begin{equation}
L\left( V^{\ast }+F_{1}^{\ast }\right) =0.  \label{7.3}
\end{equation}%
Let $\delta \in \left( 0,1\right) $ be the level of noise in the data. We
assume that%
\begin{equation}
\left\Vert F_{1}-F_{1}^{\ast }\right\Vert _{H_{N}^{2}\left( \Omega \right)
}\leq \delta .  \label{7.30}
\end{equation}%
Hence, 
\begin{equation}
\left\Vert L\left( F_{1}\right) -L\left( F_{1}^{\ast }\right) \right\Vert
_{L_{2,N}\left( \Omega \right) }\leq B\delta .  \label{7.4}
\end{equation}%
Let $\varepsilon \in \left( 0,A\right) $ be an arbitrary number. Denote 
\[
\Omega _{\varepsilon }=\left\{ \mathbf{x}=\left( x,y,z\right) :-A<x,y<A,z\in
\left( -A,A-\varepsilon \right) \right\} \subset \Omega . 
\]

\textbf{Theorem 7.2 }(convergence rate of regularized solutions).\emph{\ Let 
}$\gamma =\gamma \left( \delta \right) =\delta ^{2}$\emph{. Let }$V_{\min
,\delta }\in H_{0,bsc,N}^{3}\left( \Omega \right) $\emph{\ be the minimizer
of the functional }$J_{\delta ^{2}}\left( V\right) ,$\emph{\ which is
guaranteed by Theorem 7.1. Denote} $U_{\min ,\delta }=V_{\min ,\delta
}+F_{1}.$ \emph{There exists a sufficiently small number }$\delta _{0}\in
\left( 0,1\right) $\emph{\ and a number }$\alpha \in \left( 0,1/2\right) $%
\emph{\ such that for all }$\delta \in \left( 0,\delta _{0}\right) $\emph{\
the following convergence rate holds:}%
\begin{equation}
\left\Vert U_{\min ,\delta }-U^{\ast }\right\Vert _{H_{N}^{2}\left( \Omega
_{\varepsilon }\right) }\leq B\left( 1+\left\Vert U^{\ast }\right\Vert
_{H_{N}^{3}\left( \Omega \right) }+\left\Vert F_{1}^{\ast }\right\Vert
_{H_{N}^{3}\left( \Omega \right) }\right) \delta ^{\alpha }.  \label{7.5}
\end{equation}

\textbf{Remark 7.1}. \emph{The subdomain }$\Omega _{\varepsilon }\subset
\Omega ,$\emph{\ where convergence estimate (\ref{7.5}) is given, is larger
than a paraboloidal subdomain of }$\Omega ,$\emph{\ where the convergence
estimate of Theorem 3.2 of \cite{KQRM} is given in the elliptic case. }

\textbf{Proof of Theorem 7.2}. It follows from the variational principle that%
\begin{equation}
\left( L\left( V_{\min ,\delta }\right) ,L\left( h\right) \right) +\delta
^{2}\left[ V_{\min ,\delta },h\right] =-\left( L\left( F_{1}\right) ,L\left(
h\right) \right) ,\text{ }\forall h\in H_{0,bsc,N}^{3}\left( \Omega \right) .
\label{7.6}
\end{equation}%
Also, by (\ref{7.3})%
\begin{equation}
\left( L\left( V^{\ast }\right) ,L\left( h\right) \right) +\delta ^{2}\left[
V^{\ast },h\right] =-\left( L\left( F_{1}^{\ast }\right) ,L\left( h\right)
\right) +\delta ^{2}\left[ V^{\ast },h\right] ,\text{ }\forall h\in
H_{0,bsc,N}^{3}\left( \Omega \right) .  \label{7.7}
\end{equation}%
Denote $\widetilde{V}_{\delta }=V_{\min ,\delta }-V^{\ast }$ and $G_{\delta
}=L\left( F_{1}^{\ast }\right) -L\left( F_{1}\right) .$ By (\ref{7.4})%
\begin{equation}
\left\Vert G_{\delta }\right\Vert _{L_{2,N}\left( \Omega \right) }\leq
B\delta .  \label{7.8}
\end{equation}%
Subtract (\ref{7.7}) from (\ref{7.6}). We obtain%
\begin{equation}
\left( L\left( \widetilde{V}_{\delta }\right) ,L\left( h\right) \right)
+\delta ^{2}\left[ \widetilde{V}_{\delta },h\right]  \label{7.9}
\end{equation}%
\[
=\left( G_{\delta },L\left( h\right) \right) -\delta ^{2}\left[ V^{\ast },h%
\right] ,\text{ }\forall h\in H_{0,bsc,N}^{3}\left( \Omega \right) . 
\]%
Setting in (\ref{7.9}) $h=\widetilde{V}_{\delta },$ we obtain%
\begin{equation}
\left\Vert L\left( \widetilde{V}_{\delta }\right) \right\Vert
_{L_{2,N}\left( \Omega \right) }^{2}+\delta ^{2}\left\Vert \widetilde{V}%
_{\delta }\right\Vert _{H_{N}^{3}\left( \Omega \right) }^{2}=\left(
G_{\delta },L\left( \widetilde{V}_{\delta }\right) \right) -\delta ^{2}\left[
V^{\ast },\widetilde{V}_{\delta }\right] .  \label{7.10}
\end{equation}%
By the Cauchy-Schwarz inequality and (\ref{7.8}),%
\[
\left( G_{\delta },L\left( \widetilde{V}_{\delta }\right) \right) -\delta
^{2}\left[ V^{\ast },\widetilde{V}_{\delta }\right] \leq B\frac{\delta ^{2}}{%
2}+\frac{1}{2}\left\Vert L\left( \widetilde{V}_{\delta }\right) \right\Vert
_{L_{2,N}\left( \Omega \right) }^{2} 
\]%
\[
+\frac{\delta ^{2}}{2}\left\Vert \widetilde{V}_{\delta }\right\Vert
_{H_{N}^{3}\left( \Omega \right) }^{2}+\frac{\delta ^{2}}{2}\left\Vert
V^{\ast }\right\Vert _{H_{N}^{3}\left( \Omega \right) }^{2}. 
\]%
Comparison of this with (\ref{7.10}) gives%
\[
\left\Vert L\left( \widetilde{V}_{\delta }\right) \right\Vert
_{L_{2,N}\left( \Omega \right) }^{2}+\delta ^{2}\left\Vert \widetilde{V}%
_{\delta }\right\Vert _{H_{N}^{3}\left( \Omega \right) }^{2}\leq B\left(
1+\left\Vert V^{\ast }\right\Vert _{H_{N}^{3}\left( \Omega \right)
}^{2}\right) \delta ^{2}. 
\]%
Hence, 
\begin{equation}
\left\Vert \widetilde{V}_{\delta }\right\Vert _{H_{N}^{3}\left( \Omega
\right) }^{2}\leq B\left( 1+\left\Vert V^{\ast }\right\Vert
_{H_{N}^{3}\left( \Omega \right) }^{2}\right) ,  \label{7.11}
\end{equation}%
\begin{equation}
\left\Vert L\left( \widetilde{V}_{\delta }\right) \right\Vert
_{L_{2,N}\left( \Omega \right) }^{2}\leq B\left( 1+\left\Vert V^{\ast
}\right\Vert _{H_{N}^{3}\left( \Omega \right) }^{2}\right) \delta ^{2}.
\label{7.12}
\end{equation}

Keeping in mind (\ref{7.11}), we will now work with (\ref{7.12}),%
\[
\left\Vert L\left( \widetilde{V}_{\delta }\right) \right\Vert
_{L_{2,N}\left( \Omega \right) }^{2}=\mathop{\displaystyle \int}%
\limits_{\Omega }\left[ L\left( \widetilde{V}_{\delta }\right) \right] ^{2}d%
\mathbf{x}=\mathop{\displaystyle \int}\limits_{\Omega }\left[ L\left( 
\widetilde{V}_{\delta }\right) \right] ^{2}\mu _{1,\lambda }^{2}\left(
z\right) e^{-2\lambda \left( z-R\right) ^{2}}d\mathbf{x} 
\]%
\[
\geq e^{-2\lambda \left( A+R\right) ^{2}}\mathop{\displaystyle \int}%
\limits_{\Omega }\left[ L\left( \widetilde{V}_{\delta }\right) \right]
^{2}\mu _{1,\lambda }^{2}\left( z\right) d\mathbf{x.} 
\]%
Hence, (\ref{7.12}) implies%
\begin{equation}
\mathop{\displaystyle \int}\limits_{\Omega }\left[ L\left( \widetilde{V}%
_{\delta }\right) \right] ^{2}\mu _{1,\lambda }^{2}d\mathbf{x\leq }B\left(
1+\left\Vert V^{\ast }\right\Vert _{H_{N}^{3}\left( \Omega \right)
}^{2}\right) e^{2\lambda \left( A+R\right) ^{2}}\delta ^{2}.  \label{7.13}
\end{equation}

We now apply the Carleman estimate of Theorem 6.1 to the left hand side of (%
\ref{7.13}) for all $\lambda \geq \lambda _{0},$%
\[
\mathop{\displaystyle \int}\limits_{\Omega }\left[ L\left( \widetilde{V}%
_{\delta }\right) \right] ^{2}\mu _{1,\lambda }^{2}d\mathbf{x\geq }\frac{1}{2%
}\mathop{\displaystyle \int}\limits_{\Omega }\left( \Delta \widetilde{V}%
_{\delta }\right) ^{2}\mu _{1,\lambda }^{2}d\mathbf{x-}B\mathop{%
\displaystyle \int}\limits_{\Omega }\left( \left\vert \nabla \widetilde{V}%
_{\delta }\right\vert ^{2}+\widetilde{V}_{\delta }^{2}\right) \mu
_{1,\lambda }^{2}d\mathbf{x} 
\]%
\begin{equation}
\geq \frac{C}{\lambda }\mathop{\displaystyle \sum }\limits_{i,j=1}^{3}%
\mathop{\displaystyle \int}\limits_{\Omega }\left( \partial _{x_{i}}\partial
_{x_{j}}\widetilde{V}_{\delta }\right) ^{2}\mu _{1,\lambda }^{2}d\mathbf{x+}%
C\lambda \mathop{\displaystyle \int}\limits_{\Omega }\left( \left\vert
\nabla \widetilde{V}_{\delta }\right\vert ^{2}+\lambda ^{2}\widetilde{V}%
_{\delta }^{2}\right) \mu _{1,\lambda }^{2}d\mathbf{x}  \label{7.14}
\end{equation}%
\[
\mathbf{-}C\lambda ^{3}e^{2\lambda \left( A-R\right) ^{2}}\left\Vert 
\widetilde{V}_{\delta }\right\Vert _{H_{N}^{3}\left( \Omega \right) }^{2}-B%
\mathop{\displaystyle \int}\limits_{\Omega }\left( \left\vert \nabla 
\widetilde{V}_{\delta }\right\vert ^{2}+\widetilde{V}_{\delta }^{2}\right)
\mu _{1,\lambda }^{2}\left( z\right) d\mathbf{x.} 
\]%
Let the number $\lambda _{1}\geq \lambda _{0}>1$ be such that $C\lambda
_{1}>2B.$ Then we obtain from (\ref{7.14}) for all $\lambda \geq \lambda
_{1} $ and with a different constant $C$%
\[
\mathop{\displaystyle \int}\limits_{\Omega }\left[ L\left( \widetilde{V}%
_{\delta }\right) \right] ^{2}\mu _{1,\lambda }^{2}\left( z\right) d\mathbf{%
x\geq }\frac{C}{\lambda }\mathop{\displaystyle \sum }\limits_{i,j=1}^{3}%
\mathop{\displaystyle \int}\limits_{\Omega }\left( \partial _{x_{i}}\partial
_{x_{j}}\widetilde{V}_{\delta }\right) ^{2}\mu _{1,\lambda }^{2}d\mathbf{x+} 
\]%
\[
+C\lambda \mathop{\displaystyle \int}\limits_{\Omega }\left( \left\vert
\nabla \widetilde{V}_{\delta }\right\vert ^{2}+\lambda ^{2}\widetilde{V}%
_{\delta }^{2}\right) \mu _{1,\lambda }^{2}d\mathbf{x-}C\lambda
^{3}e^{2\lambda \left( A-R\right) ^{2}}\left\Vert \widetilde{V}_{\delta
}\right\Vert _{H_{N}^{3}\left( \Omega \right) }^{2} 
\]%
\[
\geq e^{2\lambda \left( A-\varepsilon -R\right) ^{2}}\frac{C}{\lambda }%
\mathop{\displaystyle \sum }\limits_{i,j=1}^{3}\mathop{\displaystyle \int}%
\limits_{\Omega _{\varepsilon }}\left( \partial _{x_{i}}\partial _{x_{j}}%
\widetilde{V}_{\delta }\right) ^{2}\mu _{1,\lambda }^{2}d\mathbf{x+}%
e^{2\lambda \left( A-\varepsilon -R\right) ^{2}}C\mathop{\displaystyle \int}%
\limits_{\Omega _{\varepsilon }}\left( \left\vert \nabla \widetilde{V}%
_{\delta }\right\vert ^{2}+\widetilde{V}_{\delta }^{2}\right) d\mathbf{x-} 
\]%
\[
-C\lambda ^{3}e^{2\lambda \left( A-R\right) ^{2}}\left\Vert \widetilde{V}%
_{\delta }\right\Vert _{H_{N}^{3}\left( \Omega \right) }^{2}. 
\]%
Combining this estimate with (\ref{7.11}) and (\ref{7.12}), we obtain%
\[
\left\Vert \widetilde{V}_{\delta }\right\Vert _{H^{2}\left( \Omega
_{\varepsilon }\right) }^{2}\leq B\left( 1+\left\Vert V^{\ast }\right\Vert
_{H_{N}^{2}\left( \Omega \right) }^{2}\right) \lambda ^{4}\exp \left[
2\lambda \left( \left( A-R\right) ^{2}-\left( A-\varepsilon -R\right)
^{2}\right) \right] 
\]%
\begin{equation}
+B\left( 1+\left\Vert V^{\ast }\right\Vert _{H_{N}^{3}\left( \Omega \right)
}^{2}\right) \lambda ^{4}e^{2\lambda \left( A+R\right) ^{2}}\delta ^{2}.
\label{7.15}
\end{equation}%
Choose $\lambda _{2}\geq \lambda _{1}$ such that 
\begin{equation}
\lambda ^{4}<e^{\lambda \varepsilon ^{2}},\text{ }\forall \lambda \geq
\lambda _{2}.  \label{7.16}
\end{equation}%
Then%
\begin{equation}
\lambda ^{4}\exp \left[ 2\lambda \left( \left( A-R\right) ^{2}-\left(
A-\varepsilon -R\right) ^{2}\right) \right]  \label{7.160}
\end{equation}%
\[
<\exp \left\{ -2\lambda \varepsilon \left[ 2\left( R-A\right) +\frac{%
\varepsilon }{2}\right] \right\} . 
\]%
By (\ref{7.15})-(\ref{7.160})%
\[
\left\Vert \widetilde{V}_{\delta }\right\Vert _{H^{2}\left( \Omega
_{\varepsilon }\right) }^{2}\leq B\left( 1+\left\Vert V^{\ast }\right\Vert
_{H_{N}^{3}\left( \Omega \right) }^{2}\right) \exp \left\{ -2\lambda
\varepsilon \left[ 2\left( R-A\right) +\frac{\varepsilon }{2}\right]
\right\} 
\]%
\begin{equation}
+B\left( 1+\left\Vert V^{\ast }\right\Vert _{H_{N}^{3}\left( \Omega \right)
}^{2}\right) e^{2\lambda \left[ \left( A+R\right) ^{2}+\varepsilon ^{2}/2%
\right] }\delta ^{2}.  \label{7.17}
\end{equation}%
Choose $\lambda =\lambda \left( \delta \right) $ such that 
\begin{equation}
\exp \left[ 2\lambda \left( \left( A+R\right) ^{2}+\varepsilon ^{2}/2\right) %
\right] \delta ^{2}=\delta .  \label{7.18}
\end{equation}%
Hence, 
\begin{equation}
2\lambda \left( \delta \right) =\ln \left( \delta ^{\left[ \left( A+R\right)
^{2}+\varepsilon ^{2}/2\right] ^{-1}}\right) .  \label{7.19}
\end{equation}%
Hence, by the choice (\ref{7.19}) we assume that $\delta _{0}$ is so small
that 
\[
\ln \left( \delta _{0}^{\left[ \left( A+R\right) ^{2}+\varepsilon ^{2}/2%
\right] ^{-1}}\right) >2\lambda _{2}. 
\]%
We set $\delta \in \left( 0,\delta _{0}\right) .$ Hence,%
\[
\exp \left\{ -2\lambda \left( \delta \right) \varepsilon \left[ 2\left(
R-A\right) +\frac{\varepsilon }{2}\right] \right\} =\delta ^{2\alpha }, 
\]%
\[
2\alpha =\frac{\varepsilon \left[ 2\left( R-A\right) +\varepsilon /2\right] 
}{\left( A+R\right) ^{2}+\varepsilon ^{2}/2}. 
\]%
Since $\varepsilon \in \left( 0,A\right) ,$ then, using elementary
calculations, we obtain 
\begin{equation}
2\alpha =\frac{\varepsilon \left[ 2\left( R-A\right) +\varepsilon /2\right] 
}{\left( A+R\right) ^{2}+\varepsilon ^{2}/2}<1.  \label{7.20}
\end{equation}%
Hence, (\ref{7.300}), (\ref{7.17})-(\ref{7.19}) and the triangle inequality
imply%
\begin{equation}
\left\Vert \widetilde{V}_{\delta }\right\Vert _{H^{2}\left( \Omega
_{\varepsilon }\right) }\leq B\left( 1+\left\Vert V^{\ast }\right\Vert
_{H_{N}^{3}\left( \Omega \right) }\right) \delta ^{\alpha }  \label{7.21}
\end{equation}%
\[
\leq B\left( 1+\left\Vert U^{\ast }\right\Vert _{H_{N}^{3}\left( \Omega
\right) }+\left\Vert F_{1}^{\ast }\right\Vert _{H_{N}^{3}\left( \Omega
\right) }\right) \delta ^{\alpha }. 
\]%
Now, since%
\[
\widetilde{V}_{\delta }=V_{\min ,\delta }-V^{\ast }=\left( U_{\min ,\delta
}-F_{1}\right) -\left( U^{\ast }-F_{1}^{\ast }\right) =\left( U_{\min
,\delta }-U^{\ast }\right) +\left( F_{1}^{\ast }-F_{1}\right) 
\]%
\[
=\left( U_{\min ,\delta }-U^{\ast }\right) +\left( F_{1}^{\ast
}-F_{1}\right) , 
\]%
then (\ref{7.30}) and the triangle inequality imply that 
\begin{equation}
\left\Vert \widetilde{V}_{\delta }\right\Vert _{H^{2}\left( \Omega
_{\varepsilon }\right) }\geq \left\Vert U_{\min ,\delta }-U^{\ast
}\right\Vert _{H^{2}\left( \Omega _{\varepsilon }\right) }-\delta .
\label{7.22}
\end{equation}%
Finally, since by (\ref{7.20}) $\delta <\delta ^{\alpha },$ then, using (\ref%
{7.21}) and (\ref{7.22}), we obtain (\ref{7.5}) with a new constant $B.$ $%
\square $

\subsection{Convergence rate of functions $q\left( \mathbf{x}\right) $}

\label{q}

While Theorem 7.2 establishes convergence the rate of regularized solutions
of the Minimization Problem 1, the next natural question is about the
convergence rate of corresponding target functions $q\left( \mathbf{x}%
\right) ,$ which are the approximate solutions of the system of integral
equations (\ref{2.19}), (\ref{2.20}). Recall that these functions are
reconstructed as described in Remark 5.1. Theorem 7.3 follows immediately
from Remark 5.1 and Theorem 7.2.

\textbf{Theorem 7.3}. \emph{Assume that all conditions of subsection 7.2,
described above as well as of Theorem 7.2 hold. Let} $q_{\min ,\delta
}\left( \mathbf{x}\right) $ \emph{and} $q^{\ast }\left( \mathbf{x}\right) $ 
\emph{be the functions, which are reconstructed from the vector functions }$%
U_{\min ,\delta }$\emph{\ and }$U^{\ast }$\emph{\ respectively by the
procedure described in Remark 5.1. Let the number }$\alpha \in \left(
0,1/2\right) $\emph{\ be the same as in Theorem 7.1. Then the following
convergence rate holds:}%
\[
\left\Vert q_{\min ,\delta }-q^{\ast }\right\Vert _{L_{2}\left( \Omega
_{\varepsilon }\right) }\leq B\left( 1+\left\Vert U^{\ast }\right\Vert
_{H_{N}^{3}\left( \Omega \right) }+\left\Vert F_{1}^{\ast }\right\Vert
_{H_{N}^{3}\left( \Omega \right) }\right) \delta ^{\alpha }. 
\]

\section{Numerical Studies}

\label{8}

\subsection{Numerical implementation}

\label{8.1}

\begin{figure}[tbp]
\begin{center}
\begin{tabular}{cc}
\includegraphics[width=8cm]{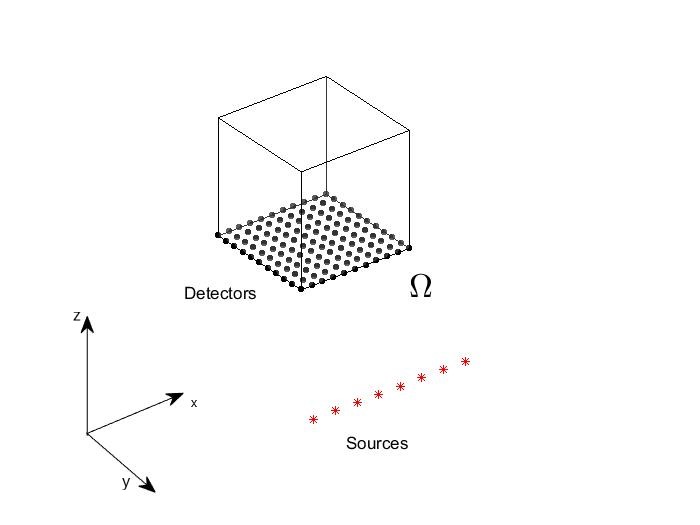} & %
\includegraphics[width=8cm]{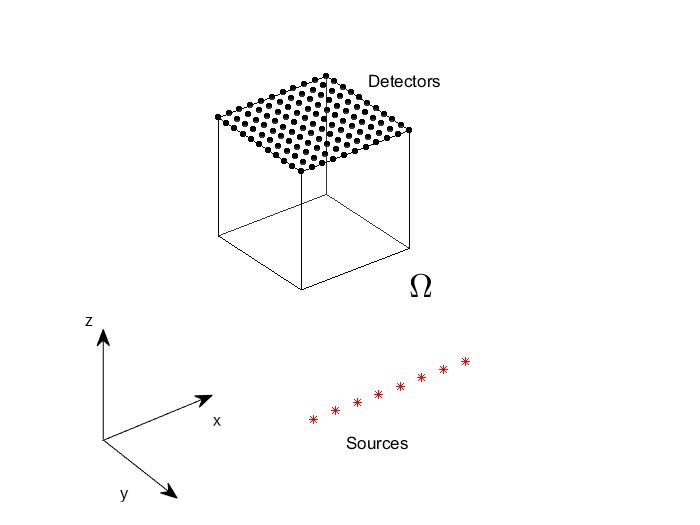} \\ 
(a)$S=S_{bsc}$ & (b) $S=S_{tr}$%
\end{tabular}%
\end{center}
\caption{\emph{A schematic diagram of domain $\Omega $, sources and
detectors.}}
\label{example0}
\end{figure}

In all the numerical tests, we have choosen the numbers $A=0.75$, $b=-5$ and 
$d=1$ in (\ref{2.0}) and (\ref{2.01}). We have used 101 sources which are
uniformly distributed on the line $L_{\text{s}}$ and $10\times 10$ detectors
uniformly distributed on the surface $S$. See Figure \ref{example0} for a
schematic diagram of domain $\Omega $, sources and detectors when $S$ is
either the backscattering part $S=S_{bsc}$ of the boundary of the cube $%
\Omega $ or its transmitted part $S=S_{tr}$.

\subsubsection{The forward problem (\protect\ref{2.3}), (\protect\ref{2.4})}

\label{8.1.1}

{To generate the data (\ref{2.14})-(\ref{2.17}) for BVP1 and BVP2 of section
4, we have solved the forward problem (\ref{2.3}), (\ref{2.4}) for each of
those 101 sources by the finite difference method. Since this procedure is
very time consuming, due to a large number of sources, we have parallelized
computations. The parallel computations were performed on 32 processors. The
delta function }$\delta \left( \mathbf{x}-\mathbf{x}_{0}\right) ${\ in (\ref%
{2.4}) was modeled as 
\[
f\left( x-x^{\left( s\right) }\right) =\left\{ 
\begin{array}{l}
\frac{1}{\varepsilon }\exp \left( -\frac{1}{1-\Vert x-x_{s}\Vert
^{2}/\varepsilon }\right) ,\text{ if }\Vert x-x_{s}\Vert ^{2}<\varepsilon ,
\\ 
0,\text{ otherwise.}%
\end{array}%
\right. 
\]%
with $\varepsilon =0.05$. In fact, we have found solution of equation (\ref%
{2.3}) with initial conditions (\ref{2.4}) in the cube} 
\[
\Phi =\left\{ \mathbf{x}=\left( x,y,z\right) :-10<x,y,z<10\right\} 
\]%
with absorbing boundary conditions on the boundary of this cube.

{The step sizes with respect to spatial and time variables were }$h=1/20${\
and }$\tau =1/100$ {respectively. We have used an implicit scheme. To
approximate integrals (\ref{2.14}), (\ref{2.15})\ over the infinite interval 
}$t\in \left( 0,\infty \right) ${\ we have continued computations of problem
(\ref{2.3}), (\ref{2.4}) for each source position }$\mathbf{x}_{0}$ {until
such a moment of time }$T${\ for which $\left\Vert u\left( \mathbf{x},%
\mathbf{x}_{0},T\right) \right\Vert_{L^\infty} \leq 10^{-3}$, where the norm is
understood in the discrete sense of grid of finite differences. We have
found that }${T=10}${\ for all locations of the point source. Next, we
replaced in integrals (\ref{2.14}), (\ref{2.15}) the integration interval }$%
t\in \left( 0,\infty \right) $ {with the interval} $t\in \left( 0,T\right) .$

\subsubsection{The inverse problem}

\label{8.1.2}

Even though we use the $H_{N}^{3}\left( \Omega \right) -$norm in the
regularization terms in (\ref{5.200}), (\ref{5.201}), we have discovered
numerically that the simpler $H_{N}^{2}\left( \Omega \right) -$norm is
sufficient. To minimize either of functionals (\ref{5.200}), (\ref{5.201}),
we have written differential operators in them in the form of finite
differences. Then we have minimized with respect to the values of the
corresponding vector functions at grid points. The grid step size was $%
\widetilde{h}=1/10$. The discrete problems of (\ref{5.200}),
(\ref{5.201}) will lead to 2 linear systems. Then we have solved the minimization problems (\ref{5.200}),
(\ref{5.201}) directly by solving these 2 linear equations. 

\subsection{Results}

\label{8.2}

In the tests of this subsection, we demonstrate the efficiency of our
numerical method. We display on Figures 2-9 the function $q\left( \mathbf{x}%
\right) =c\left( \mathbf{x}\right) -1,$ see (\ref{2.10})\textbf{. }In all
tests the background value of $q_{bkgr}\left( \mathbf{x}\right) =0$, which
corresponds to the background value of the dielectric constant $%
c_{bkgr}\left( \mathbf{x}\right) =1.$ In all tests slices are depicted to
demonstrate the values of the computed function $q\left( \mathbf{x}\right) $%
. In our computations, we took the number of Fourier coefficients $N=6$ in (%
\ref{3.8}), (\ref{3.9}).\ \ \ 

In tests 1-4, we work with the case when the surface $S$ is%
\[
S=S_{bsc}=\left\{ \mathbf{x}:\left\vert x\right\vert ,\left\vert
y\right\vert <0.75,z=-0.75\right\} 
\]%
is the backscattering part of the boundary of the cube $\Omega $, see Figure
1(a) and (\ref{4.1}). In tests 1-3, we have noiseless data on $S,$ and we
set in these tests the regularization parameter $\gamma =0$ in (\ref{5.200}%
). In test 4 we use noisy data with the noise level $\delta =0.05$ (i.e.
5\%), and we set in test 4 $\gamma =10^{-8}$ in (\ref{5.200}).

\begin{figure}[tbp]
\begin{center}
\begin{tabular}{cc}
\includegraphics[width=4.8cm]{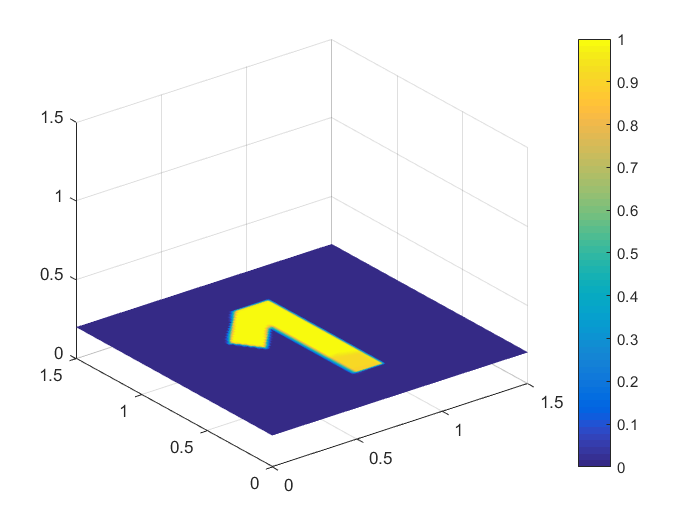} & %
\includegraphics[width=4.8cm]{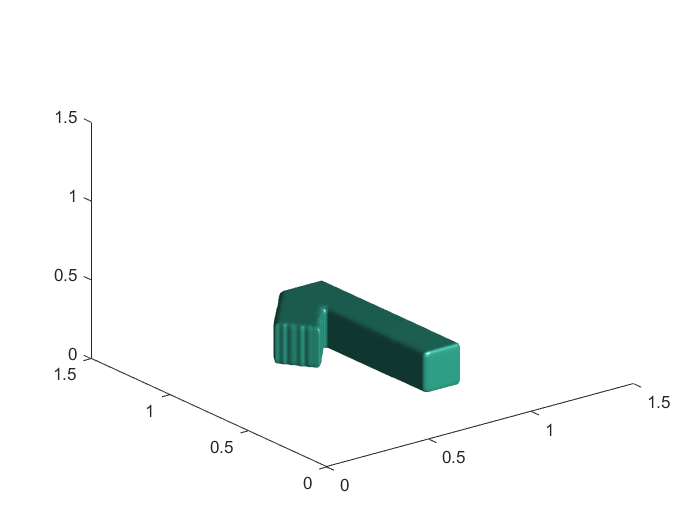} \\ 
(a)Slice image of true $q$ & (b) 3D image of true $q$ \\ 
\includegraphics[width=4.8cm]{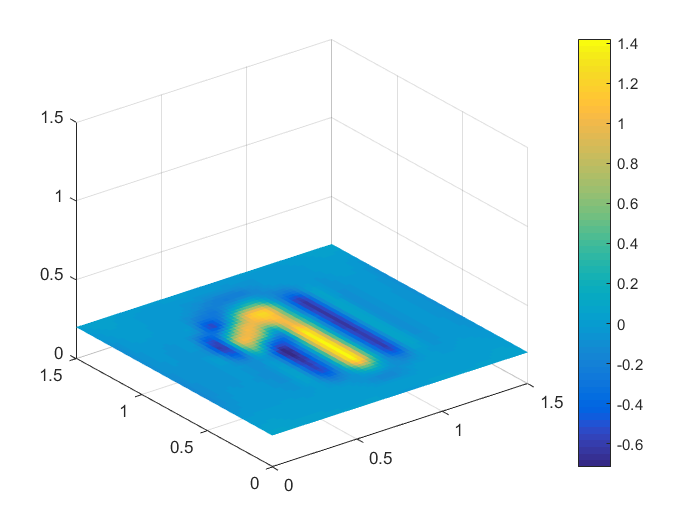} & %
\includegraphics[width=4.8cm]{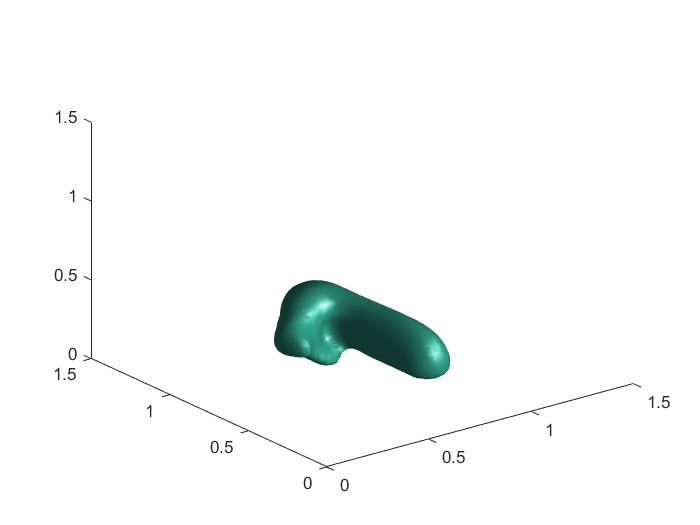} \\ 
(c) Slice image of recovered $q$ & (d) 3D image of recovered $q$%
\end{tabular}%
\end{center}
\caption{\emph{Results of Test 1. Backscattering data, $S = S_{bsc}.$
Imaging of '1' shaped $q$ with } $q =1$\emph{\ in it and $q =0$ outside. (a)
and (b) Correct images. (c) and (d) Computed images.}}
\label{example1}
\end{figure}

\textbf{Test 1}. First, we test the reconstruction by our method of the case
when the shape of our inclusion is the same as the shape of the number `1'.
The function $q(\mathbf{x})$ is depicted on Figures \ref{example1} (a) and
(b). $q=1$ inside of this inclusion and $q=0$ outside of it. See Figures \ref%
{example1} for the reconstruction results.

\begin{figure}[tbp]
\begin{center}
\begin{tabular}{cc}
\includegraphics[width=4.8cm]{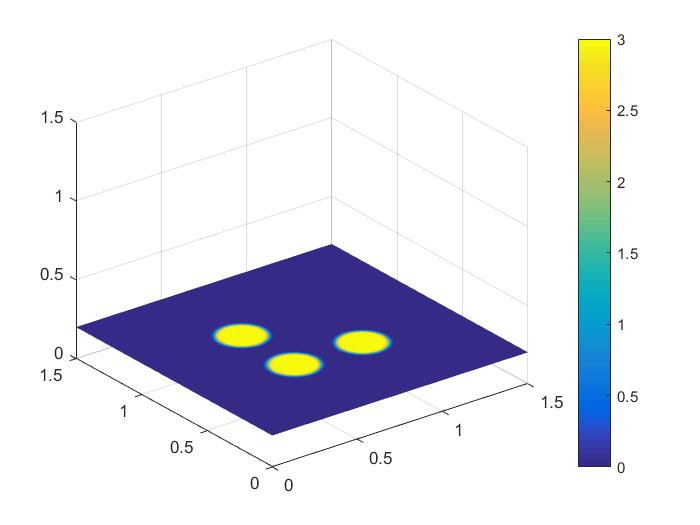} & %
\includegraphics[width=4.8cm]{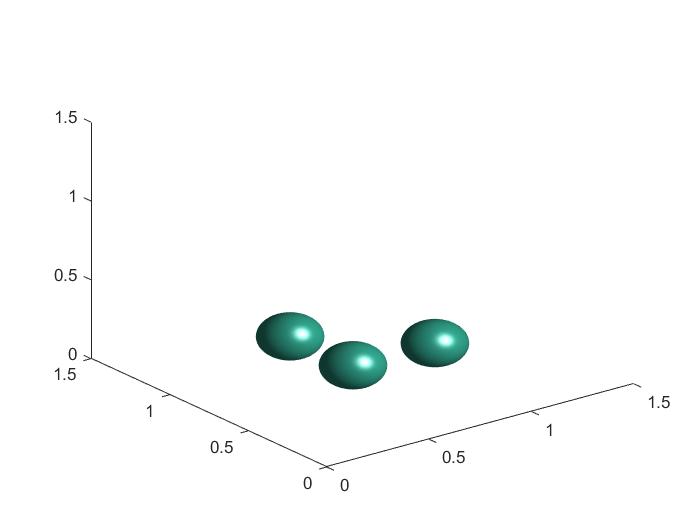} \\ 
(a) Slice image of true $q$ & (b) 3D image of true $q$ \\ 
\includegraphics[width=4.8cm]{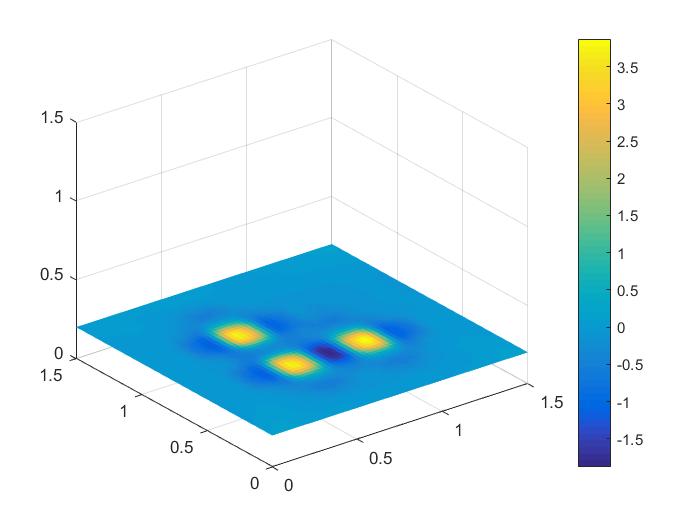} & %
\includegraphics[width=4.8cm]{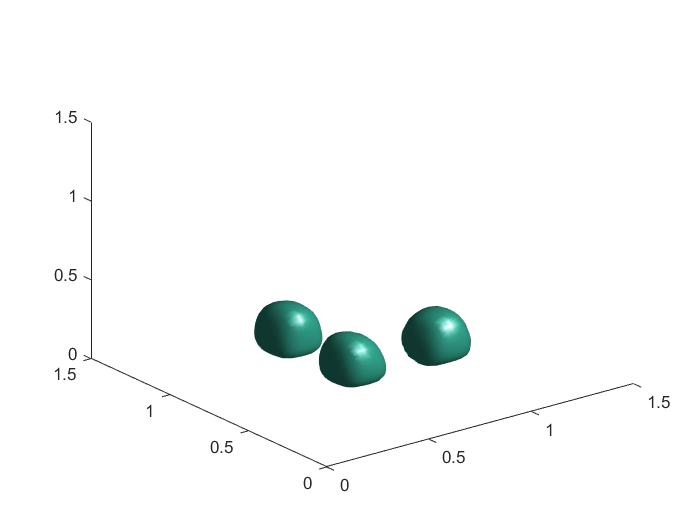} \\ 
(c) Slice image of recovered $q$ & (d) 3D image of recovered $q$%
\end{tabular}%
\end{center}
\caption{\emph{Results of Test 2. Backscattering data, $S = S_{bsc}.$\ Imaging of three ball shaped inclusions
with } $q =3$\emph{\ in each of them and $q =0$ outside. a) and b) Correct
images. c) and d) Computed images.}}
\label{example2}
\end{figure}

\textbf{Test 2}. We test the reconstruction by our method of three (3)
single ball shaped inclusions with a high contrast, see Figures \ref%
{example2} (a) and (b). $q=3$ inside of this inclusions and $q=0$ outside of
them. See Figures \ref{example2} for the results of the reconstruction.

\begin{figure}[tbp]
\begin{center}
\begin{tabular}{cc}
\includegraphics[width=4.8cm]{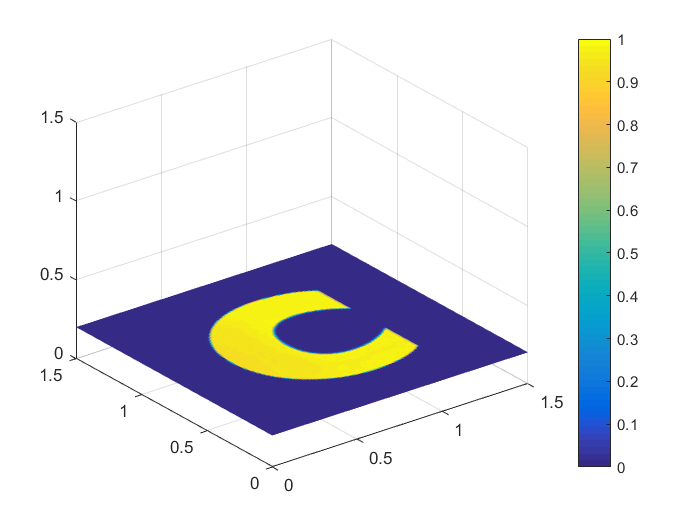} & %
\includegraphics[width=4.8cm]{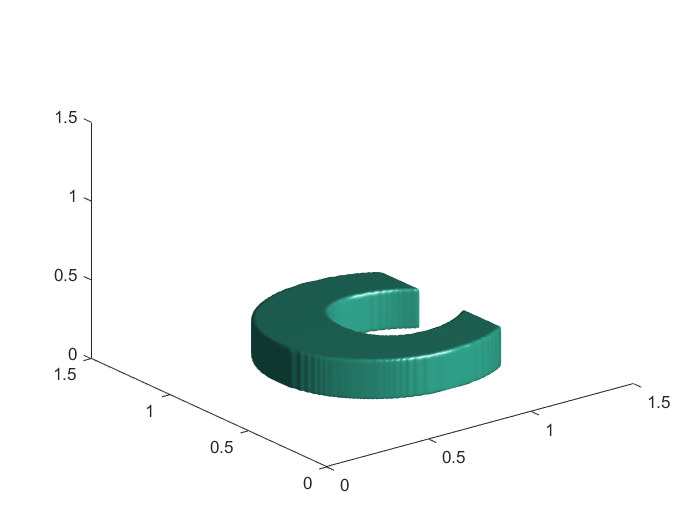} \\ 
(a) Slice image of true $q$ & (b) 3D image of true $q$ \\ 
\includegraphics[width=4.8cm]{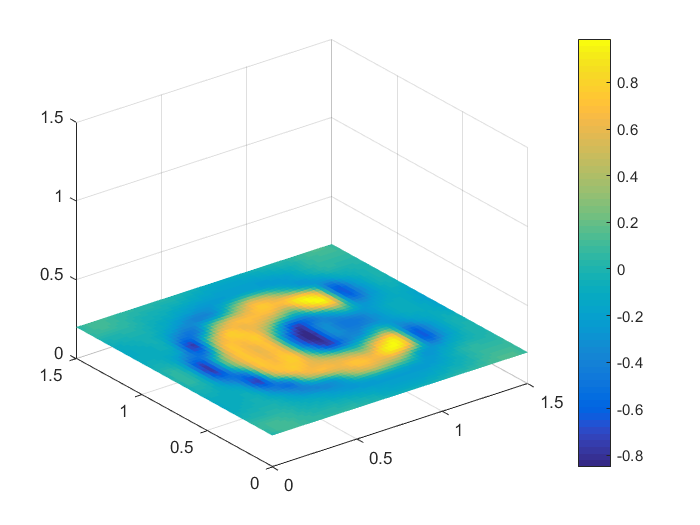} & %
\includegraphics[width=4.8cm]{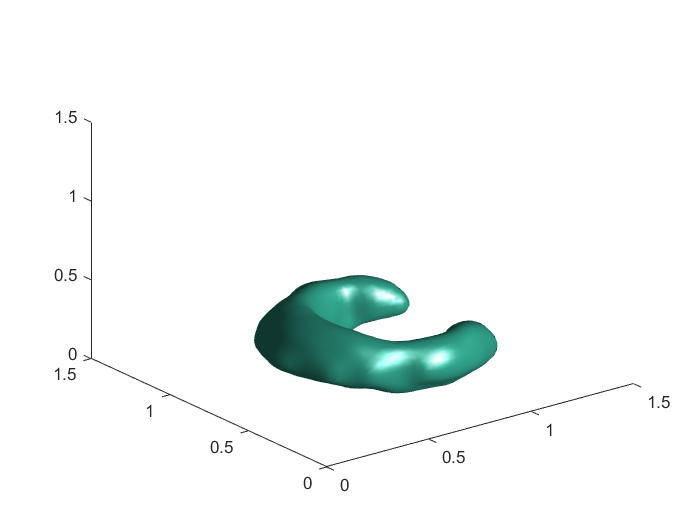} \\ 
(c) Slice image of recovered $q$ & (d) 3D image of recovered $q$%
\end{tabular}%
\end{center}
\caption{\emph{Results of Test 3. Backscattering data, $S = S_{bsc}.$ \
Imaging of 'C' shaped $q$ with } $q =1$\emph{\ in it and $q =0$ outside. (a)
and (b) Correct images. (c) and (d) Computed images.}}
\label{example3}
\end{figure}

\textbf{Test 3}. We test the reconstruction by our method of the case when
the shape of our inclusion is the same as the shape of the letter `C'. The
function $q(\mathbf{x})$ is depicted on Figures \ref{example3} (a) and (b). $%
q=1$ inside of this inclusion and $q=0$ outside of it. See Figures \ref%
{example3} for the reconstruction results.

\textbf{Test 4}. In this example, we test the stability of our algorithm
with respect to the random noise in the data. We use the same function $%
q\left( \mathbf{x}\right) $ as the one in Test 3. The random noise is added
on $S_{bsc}$ as: 
\[
g_{0,noise}\left( \mathbf{x},\mathbf{x}_{0}\right) =g_{0}\left( \mathbf{x},%
\mathbf{x}_{0}\right) +\delta \max_{\mathbf{x}\in S_{bsc},\mathbf{x}_{0}\in
L_{\text{s}}}|g_{0}\left( \mathbf{x},\mathbf{x}_{0}\right) |\xi _{{x}_{0}},%
\text{ }\mathbf{x}\in S_{bsc},\mathbf{x}_{0}\in L_{\text{s}}, 
\]%
\[
g_{1,noise}\left( \mathbf{x},\mathbf{x}_{0}\right) =g_{1}\left( \mathbf{x},%
\mathbf{x}_{0}\right) +\delta \max_{\mathbf{x}\in S_{bsc},\mathbf{x}_{0}\in
L_{\text{s}}}|g_{1}\left( \mathbf{x},\mathbf{x}_{0}\right) |\xi _{{x}_{0}},%
\text{ }\mathbf{x}\in S_{bsc},\mathbf{x}_{0}\in L_{\text{s}}, 
\]%
where functions $g_{0}(\mathbf{x},\mathbf{x}_{0}),g_{1}(\mathbf{x},\mathbf{x}%
_{0})$ are defined in (\ref{2.16}), (\ref{2.17}), $\delta $ is the noise
level and $\xi _{\mathbf{x}_{0}}$ is a random variable depending only on the
source position $\mathbf{x}_{0}$ and uniformly distributed on $[-1,1]$. We
took $\delta =0.05$, which corresponds to 5\% noise level. See Figures \ref%
{example4} for the reconstruction results.

In tests 5-8, we show the reconstruction results for the case when the
surface $S$ is 
\[
S=S_{tr}=\left\{ \mathbf{x}:\left\vert x\right\vert ,\left\vert y\right\vert
<A,z=A\right\} 
\]
is the transmitted part of the boundary of the cube $\Omega $, see Figure
1(b) and (\ref{4.2}). In tests 5-7, the we use noiseless data $S$ are set
the regularization parameter $\gamma =0$ in (\ref{5.201}). In test 8 we
consider the noisy data in the same manner as in test 5 with $\delta =0.05$
and we set in (\ref{5.201}) $\gamma =10^{-8}$.

\textbf{Test 5}. First, we test the reconstruction by our method of the case
when the shape of our inclusion is the same as the shape of the number `1'.
The function $q(\mathbf{x})$ is depicted on Figures \ref{example5} (a) and
(b). $q=1$ inside of `1' and $q=0$ outside of it. See Figures \ref{example5}
for the reconstruction results.

\textbf{Test 6}. We test the reconstruction by our method of three (3) ball
shaped inclusions with high contrast, see Figures \ref{example6} (a), (b). $%
q=3$ inside of these inclusions and $q=0$ outside of them. See Figures \ref%
{example6} for the reconstruction results.

\textbf{Test 7}. We test the reconstruction by our method of the case when
the shape of our inclusion is the same as the shape of the letter `C'. The
function $q(\mathbf{x})$ is depicted on Figures \ref{example7} (a) and (b). $%
q=1$ inside of `C' and $q=0$ outside. See Figures \ref{example7} for the
reconstruction results.

\textbf{Test 8}. In this test, we act the same way as in test 4, except that
we have replaced the surface $S=S_{bsc}$ with $S=S_{tr}$. The function $q$
is depicted in Figures \ref{example8} (a) and (b). We took the noise level $%
\delta =0.05,$ i.e. $5\%$. Figures \ref{example8} display the reconstruction
results.

\begin{figure}[tbp]
\begin{center}
\begin{tabular}{cc}
\includegraphics[width=4.8cm]{Figures/letter-C-down-slice-true.png} & %
\includegraphics[width=4.8cm]{Figures/letter-C-down-surface-true.png} \\ 
(a) Slice image of true $q$ & (b) 3D image of true $q$ \\ 
\includegraphics[width=4.8cm]{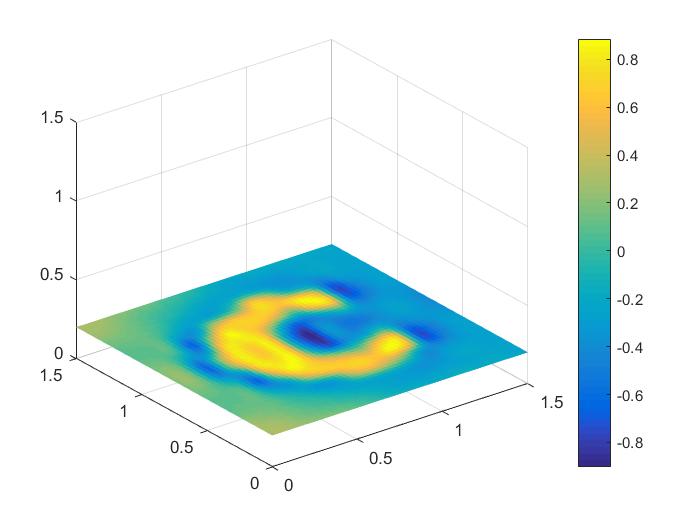}
& %
\includegraphics[width=4.8cm]{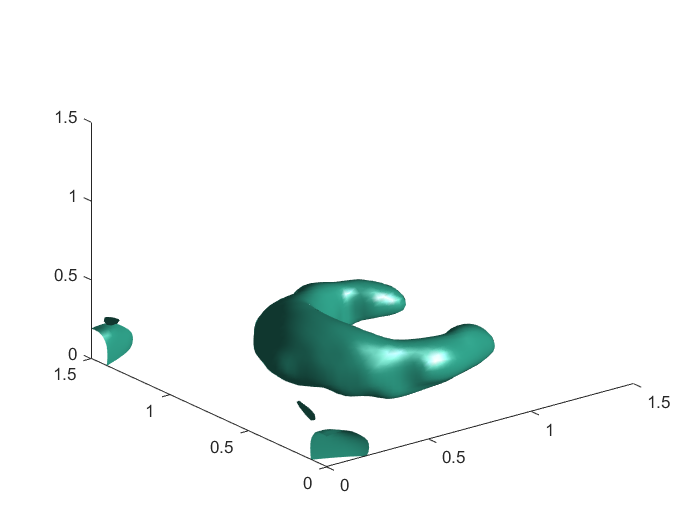}
\\ 
(c) Slice image of recovered $q$ & (d) 3D image of recovered $q$%
\end{tabular}%
\end{center}
\caption{\emph{Results of Test 4. Backscattering data, $S = S_{bsc}.$ \
Imaging of 'C' shaped $q$ with } $q =1$\emph{\ in it and $q =0$ outside. (a)
and (b) Correct images. (c) and (d) Computed images. In this case, noise
level $\protect\delta=0.05$ (i.e.5$\%$).}}
\label{example4}
\end{figure}

\begin{figure}[tbp]
\begin{center}
\begin{tabular}{cc}
\includegraphics[width=4.8cm]{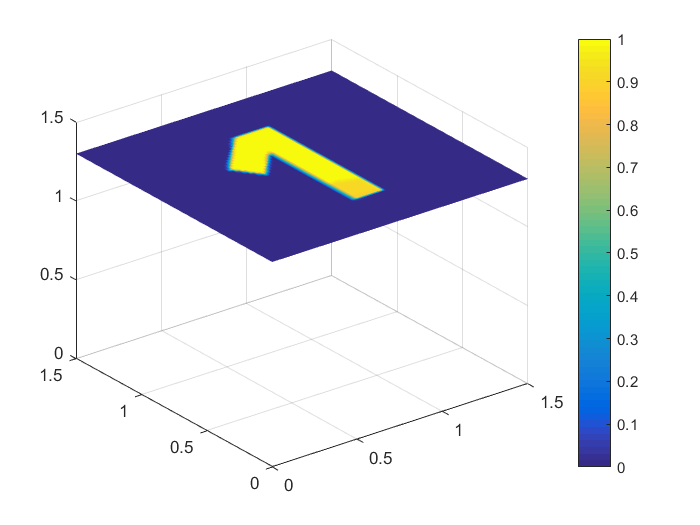} & %
\includegraphics[width=4.8cm]{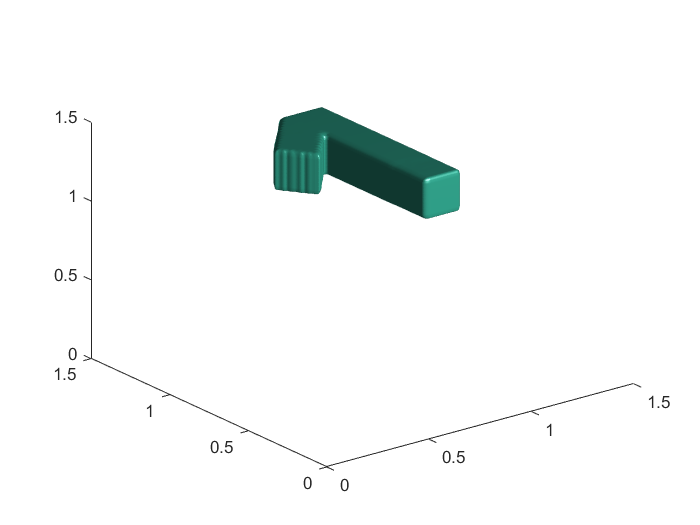} \\ 
(a)Slice image of true $q$ & (b) 3D image of true $q$ \\ 
\includegraphics[width=4.8cm]{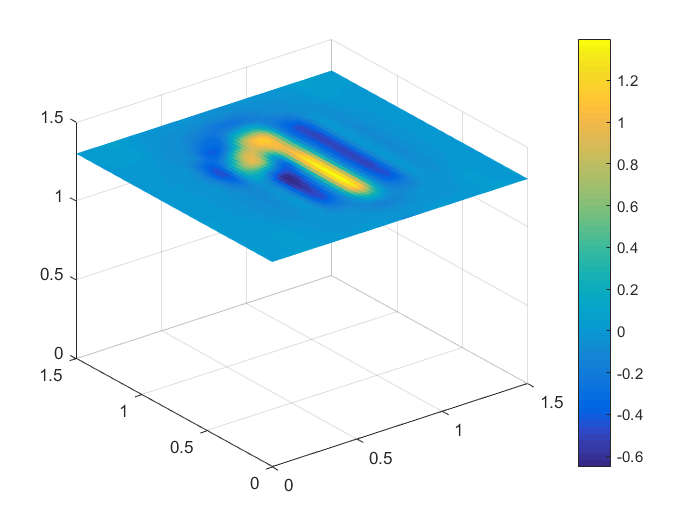} & %
\includegraphics[width=4.8cm]{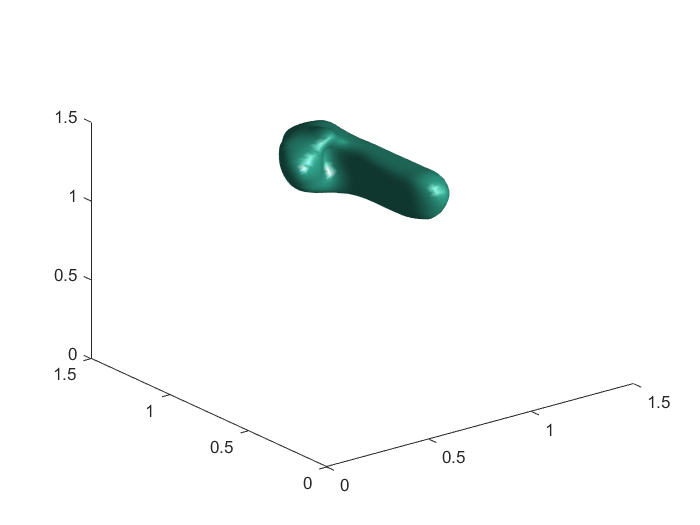} \\ 
(c) Slice image of recovered $q$ & (d) 3D image of recovered $q$%
\end{tabular}%
\end{center}
\caption{\emph{Results of Test 5. Transmitted data, $S = S_{tr}$.\ Imaging
of '1' shaped $q$ with } $q =1$\emph{\ in it and $q=0$ outside. (a) and (b)
Correct images. (c) and (d) Computed images.}}
\label{example5}
\end{figure}

\begin{figure}[tbp]
\begin{center}
\begin{tabular}{cc}
\includegraphics[width=4.8cm]{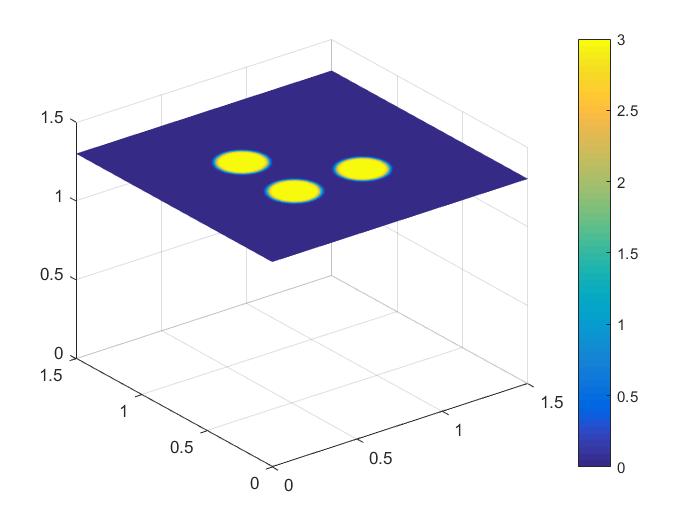} & %
\includegraphics[width=4.8cm]{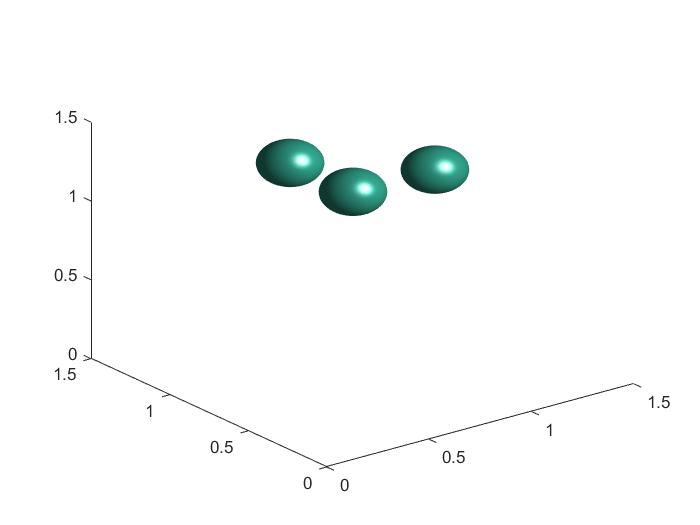} \\ 
(a) Slice image of true $q$ & (b) 3D image of true $q$ \\ 
\includegraphics[width=4.8cm]{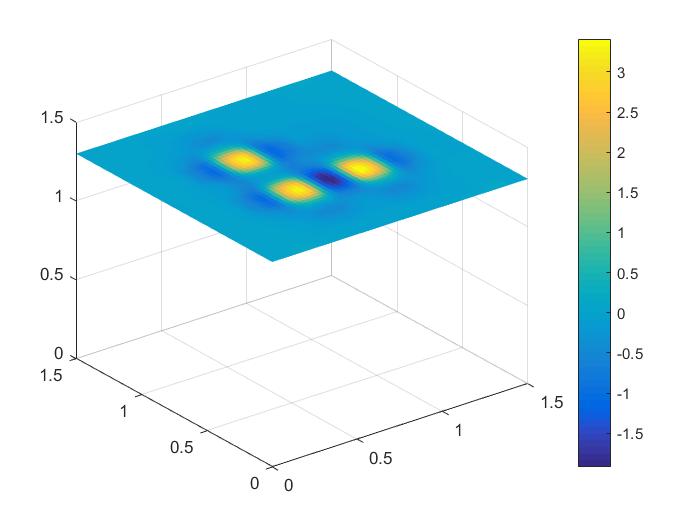} & %
\includegraphics[width=4.8cm]{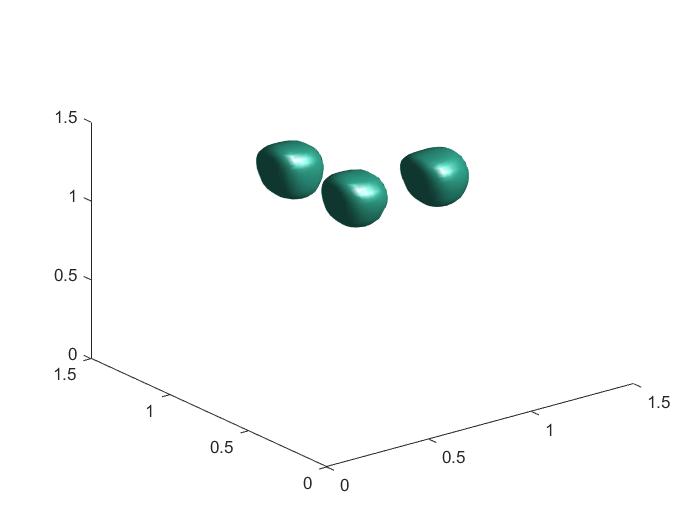} \\ 
(c) Slice image of recovered $q$ & (d) 3D image of recovered $q$%
\end{tabular}%
\end{center}
\caption{\emph{Results of Test 6. Transmitted data, $S = S_{tr}$.\ Imaging of 3 ball shaped inclusions with }
$q =3$\emph{\ in each of them and $q =0$ outside. a) and b) Correct images.
c) and d) Computed images.}}
\label{example6}
\end{figure}

\begin{figure}[tbp]
\begin{center}
\begin{tabular}{cc}
\includegraphics[width=4.8cm]{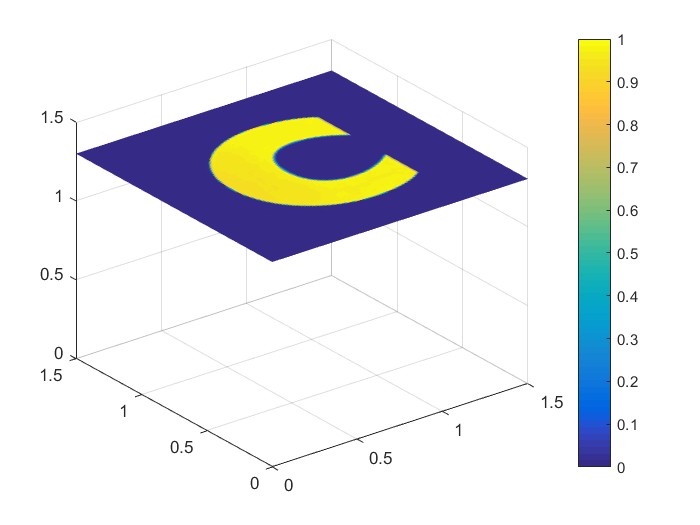} & %
\includegraphics[width=4.8cm]{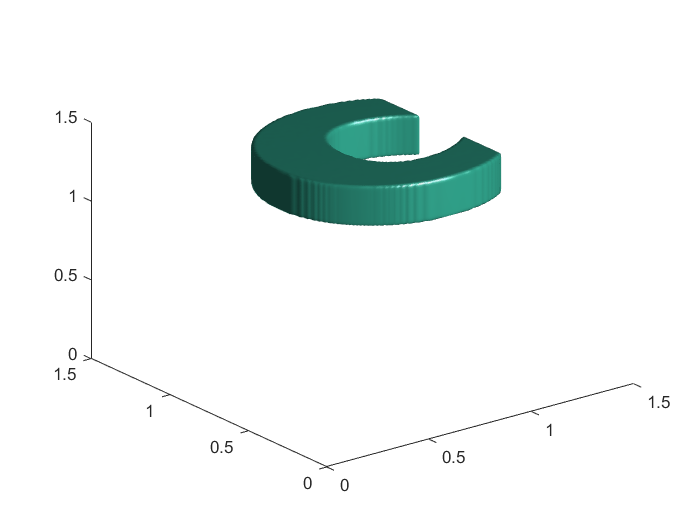} \\ 
(a) Slice image of true $q$ & (b) 3D image of true $q$ \\ 
\includegraphics[width=4.8cm]{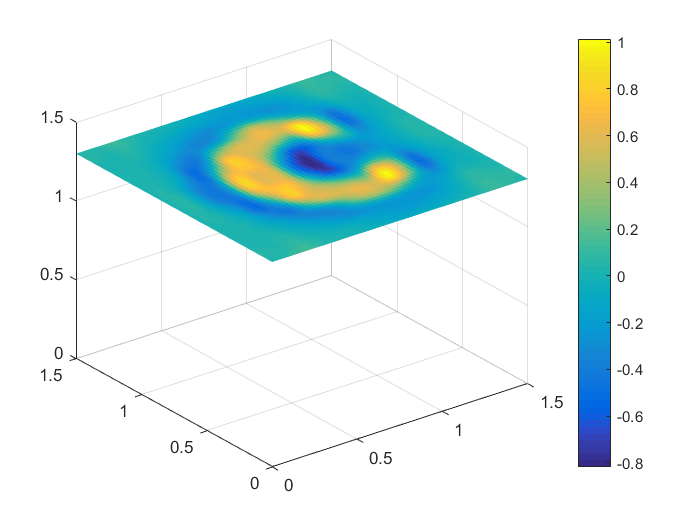} & %
\includegraphics[width=4.8cm]{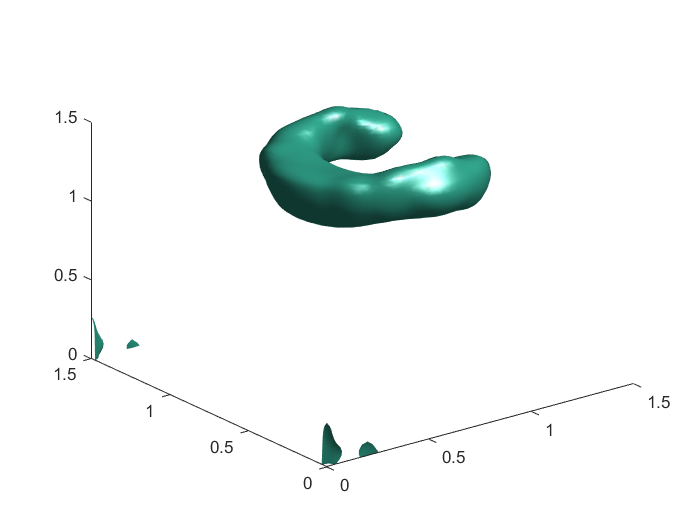} \\ 
(c) Slice image of recovered $q$ & (d) 3D image of recovered $q$%
\end{tabular}%
\end{center}
\caption{\emph{Results of Test 7. Transmitted data, $S = S_{tr}$.\ Imaging
of 'C' shaped $q$ with } $q =1$\emph{\ in it and $q=0$ outside. (a) and (b)
Correct images. (c) and (d) Computed images.}}
\label{example7}
\end{figure}

\begin{figure}[tbp]
\begin{center}
\begin{tabular}{cc}
\includegraphics[width=4.8cm]{Figures/letter-C-up-slice-true.png} & %
\includegraphics[width=4.8cm]{Figures/letter-C-up-surface-true.png} \\ 
(a) Slice image of true $q$ & (b) 3D image of true $q$ \\ 
\includegraphics[width=4.8cm]{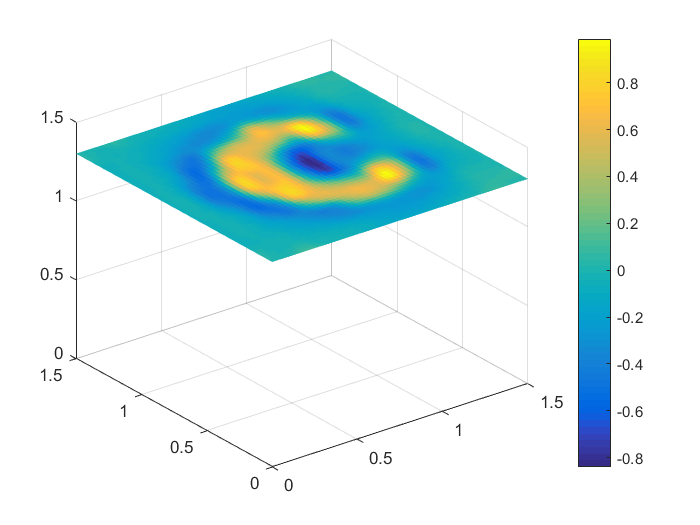}
& %
\includegraphics[width=4.8cm]{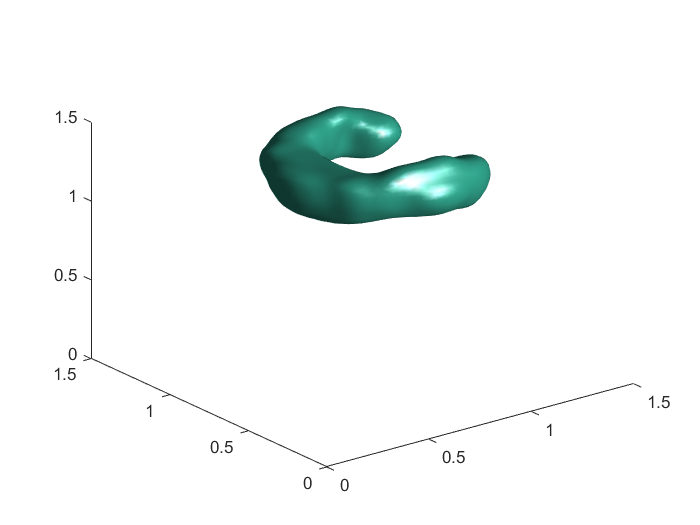}
\\ 
(c) Slice image of recovered $q$ & (d) 3D image of recovered $q$%
\end{tabular}%
\end{center}
\caption{\emph{Results of Test 8. Transmitted data, $S = S_{tr}$.\ Imaging
of 'C' shaped $q$ with } $q =1$\emph{\ in it and $q=0$ outside. (a) and (b)
Correct images. (c) and (d) Computed images. In this case, noise level $%
\protect\delta=0.05$ (i.e.5$\%$).}}
\label{example8}
\end{figure}

\end{document}